%% file: main.tex
\documentclass[11pt,a4paper]{article}
\usepackage[T1]{fontenc}

\usepackage[a4paper,top=2cm,bottom=2cm,left=3cm,right=3cm,marginparwidth=1.75cm]{geometry}
 
\usepackage{amsmath}
\usepackage{graphicx}
\usepackage{hyperref}
\usepackage{eurosym}
\usepackage{amssymb}
 
\usepackage{eurosym}
\usepackage{textcomp}
\usepackage{ulem}
\usepackage{subcaption}
\usepackage[utf8]{inputenc}
\usepackage{pgfplots}
\usepackage{xurl}
\usepackage{bm}
\usepackage{bbm}

\usepackage{apacite}
\usepackage{natbib}

\usepackage[affil-it]{authblk}

\makeatletter
\renewcommand{\section}{\@startsection {section}{1}{\z@}%
              {24pt}{12pt} {\large\scshape\bfseries}}

\renewcommand{\subsection}{\@startsection {subsection}{2}{\z@}%
             {12pt}{12pt}  {\itshape\bfseries}}

\title{\bfseries \normalsize Backcasting the Optimal Decisions in Transport Systems: An Example with Electric Vehicle Purchase Incentives}

\author{Dr. Vinith Lakshmanan*, Xavier Guichet, and Antonio Sciarretta}

\affil{IFP Energies nouvelles, France}

\date{\vspace{-5ex}}

\begin{document}
\maketitle

\section*{Short summary}\small
This study represents a first attempt to build a backcasting methodology to identify the optimal policy roadmaps in transport systems. In this methodology, desired objectives are set by decision makers at a given time horizon, and then the optimal combinations of policies to achieve these objectives are computed as a function of time (i.e., ``backcasted''). This approach is illustrated on the transportation sector by considering a specific subsystem with a single policy decision. The  subsystem describes the evolution of the passenger car fleet within a given region and its impact on greenhouse gas emissions. The optimized policy is a monetary incentive for the purchase of electric vehicles while minimizing the total budget of the state and achieving a desired CO$_2$ target. A case study applied to Metropolitan France is presented to illustrate the approach. Additionally, alternative policy scenarios are also analyzed to provide further insights.


\textbf{Keywords}: Backcasting, Transport systems, Policy making, Environmental Impacts

\section{Introduction}

The European Union's (EU) goal of carbon neutrality by 2050 requires a 90\% reduction in emissions from the transport sector compared to 1990 levels. 
The European Commission has therefore adopted a package of proposals to put EU policies on track to reduce net greenhouse gas emissions by at least 55\% by 2030 (\cite{GreenDeal}).

Governance, policies and incentives (``decisions'') play an important role in shaping the transport systems of the future, influencing the development and implementation of different technologies and modes of transport. It is therefore important to examine how decisions can best be used to steer them in the desired direction of decarbonisation. 


The traditional approach to finding policy roadmaps involves designing prospective scenarios, simulating their impacts, and drawing conclusions about effective decisions. However, this method limits choices to predefined scenarios, which represent only a small subset of possibilities. As a result, the optimum may not be reached since only selected scenarios are evaluated.

To overcome these limitations, a novel (``backcasting'') paradigm is supported in this work. In this approach, desired objectives are set by decision makers at a given time horizon, and the optimal combinations of policies to achieve these objectives are computed as a function of time (``backcasted''). 
In this way, the a priori choice of scenarios is replaced by a fully dynamic optimisation process that can explore all possible combinations.

The backcasting paradigm has been introduced since the last century (\cite{Robinson, Bibri}). It has mainly been used in a qualitative way (\cite{Levitate}) or in a quantitative way with a static optimisation procedure (\cite{Gomi, Ashina}).
However, this process can be more effectively cast as an optimal control problem, with a suitable definition of an objective function, a horizon, local and terminal constraints, etc.

Of course, since future impacts must be predicted, the new backcasting paradigm is still based on a simulation model. 
 
This model must be able to describe transport as a system, with controllable inputs, exogenous inputs, outputs, and states. The manipulable inputs represent the decisions to be optimized, which may concern local and state authorities, the EU, or even private companies. The exogenous inputs represent the influence of other related systems, such as energy, urban, economic, demographic.
The outputs represent the desired effects or the constraints to be imposed on the backcasting process. Finally, the states are the dynamics associated with the internal variables. 

In this paper, we illustrate the backcasting paradigm applied to the transport sector by considering a specific subsystem with a single decision variable. The considered subsystem describes the evolution of the passenger car fleet within a given region and its impact on greenhouse gas emissions. The optimized decision is a monetary incentive for the purchase of electric vehicles.

The prediction of fleet composition is the subject of a large body of literature (\cite{TREMOVE, High-Tool, ITF, ADEME}). Typically, dynamic fleet models are based on the evaluation of stocks and sales of different types of vehicles per time period. Stocks change over time due to the disposal of old vehicles (due to scrapping, export, change of use, etc.) and sales of new vehicles. The latter, in turn, are induced by transport demand (vkm) and mileage, and are split among vehicle types using discrete choice models (\cite{train, benakiva}). 

The GHG emissions of a given vehicle fleet are typically assessed using emission factors. CO$_2$ emissions of light-duty vehicles are regulated in the EU. Similar regulations are about to be applied to heavy-duty vehicles.

Recent studies that include electric vehicles have applied a fleet model to predict the future transportation emissions in France (\cite{ADEME, ITF}), Norway (\cite{thorne}), Japan (\cite{sato}), and the U.S. (\cite{woody}).


The paper is organized as follows. Section~\ref{sec:meth} outlines the backcasting methodology, including system model assumptions, a reduced-order model of passenger car fleet evolution, and the backcasting problem formulation with its semi-analytical solution. A full-order model and its OCP formulation are then presented. Section~\ref{sec:results} discusses case study applied to France targeting 2050 CO\(_2\) emissions using the full-order model. Finally, Sect.~\ref{sec:conclusion} provides conclusions and suggests directions for making it more realistic.

\section{Methodology}\label{sec:meth}

\subsection{Assumptions}

In this study, the backcasting process optimizes the government's monetary incentive $u(t)$ (\euro) for the purchase of electric vehicles (EV) to meet a target for greenhouse gas emissions (CO$_2$) in year $T$. The system is represented by an aggregate dynamic model that computes CO$_2$ emissions over time $t$ with a one-year time step. 

We consider only a single zone of interest and only private vehicles as the mode of transportation. The latter is composed of two vehicle types (thermal, $v=1$ and electric, $v=2$). Additionally, the vehicles are distinguished by $A+1$ age classes ($a=0\dots A$), where $a=0$ corresponds to new vehicles sold. The vehicle stock by type and age in year $t$ is denoted as $S_{va}(t)$, and its survival rate over time due to factors such as scrappage, wear, and obsolescence, is modeled using an age-dependent factor $\eta_a$.

Several model inputs are considered exogenous such as transport demand $G(t)$ (vehicle-km), mileage $M(t)$ (km/y), and purchase and operating costs (i.e., fuel, maintenance, and insurance) by vehicle type (\euro). The latter are considered as the main  determinants for the choice of new vehicles, alongside the development rate of refueling infrastructure. Conversely, socio-economic factors like age, gender, and income are not explicitly considered; instead, an adoption coefficient based on the Bass model (\cite{Bass}) is used to model new technology adoption, such as the EV(\cite{macmanus, sterman}).  

Tailpipe CO$_2$ emissions $E(t)$ are described using emission factors $\epsilon_{1a}$ (g/km) for thermal vehicles distinguished by age. Obviously, $\epsilon_{2a}\equiv 0$. Additionally, cumulative emissions $\mathcal{E}(t)$ represent the sum of the yearly CO$_2$ emissions from the present until year $t$.

\subsection{Modeling and backcasting: Reduced-order model}

We first derive a simple model with constant mileage $M$, two age classes, i.e., new and old vehicles ($A=1$), and a constant survival rate $\eta$. Therefore, we remove the age index and denote the stock of new vehicle by $N(t)$ and that of old vehicles by $O(t)$.
Under these assumptions, the total stock always equals the transport demand, i.e., $S(t)=N(t)+O(t)=S_1(t)+S_2(t)\equiv G(t)/M$. 

The model can be summarized as follows
\begin{equation} \label{eqn:summ}
   E(t) = \epsilon_1 M S_1(t)\;,
\end{equation}
\begin{equation} \label{eqn:S1}
  S_1(t) = \eta S_1(t-1) + N(t)P_1(t,u(t)) \;,
\end{equation}
\begin{equation}\label{cumE}
    \mathcal{E}(t) = \mathcal{E}(t-1) + \epsilon_1 MS_1(t)\;, 
\end{equation}
\begin{equation} \label{eqn:P}
    P_1(t,u(t)) = \dfrac{\mathcal{P}(t)}{\mathcal{P}(t)+\mathcal{Q}(t) \mathcal{R}(t)^{u(t)}}\;,
\end{equation}
where
\begin{equation}
  N(t) = \dfrac{G(t)-\eta G(t-1)}{M}  
\end{equation}
\begin{equation} \label{eqn:p}
  \mathcal{P}(t)= e^{\mu U_1(t)}  
\end{equation}
\begin{equation}\label{eqn:q}
    \mathcal{Q}(t) = \exp\left(\mu(1-c_2^A(t))\left(p^P\frac{C^P_2(t)}{\overline{C}^P(t)} + p^O \frac{C^O_2(t)}{\overline{C}^O(t)} + p^I (1-c^I_2(t)) \right) \right),
\end{equation}
and 
\begin{equation} \label{eqn:R}
  \mathcal{R}(t)=\exp\left(-\mu(1-c_2^A(t))\frac{p^P}{\overline{C}^P(t)} \right)
\end{equation}
are explicit functions of time only.

The latter derive from a logit discrete choice model based on a utility function defined as
\begin{equation} \label{eqn:U}
  U_v(t) = \left(1-c_v^A(t) \right)\left(p^P\frac{C_v^P(t)-u_v(t)}{\overline{C}^P(t)} + p^O \frac{C_v^O(t)}{\overline{C}^O(t)}+p^I {(1-c_v^I(t))} \right),  
\end{equation}
%
i.e., the product of an adoption-based factor and a cost-based factor. In the latter, $C_v^P$ is the purchase price, $C^O_v$ is the sum of operating costs, $c^I_v$ is the rate of development of the refueling infrastructure (normalized to unity, by definition $c_1^I\equiv 1$), and $p$ are tuning coefficients. The average costs between the two vehicle types are given by $\overline{C}_v^P(t)$ and $\overline{C}_v^O(t)$. 
In the adoption-based utility factor, $\mu$ is a tuning coefficient and $c^A_2$ ($c^A_1\equiv 0$) is the adoption coefficient. Since the cost-based utility is negative (the coefficients $p$ are so), a prefactor less than unity increases the utility of EVs proportionally to their rate of exposure.
Since the monetary incentive applies only for EVs, we define $u_1(t) \equiv 0$, while $u_2(t) \equiv u(t)$.

By integrating (\ref{eqn:S1}) over time, the evolution of vehicle stock at time $T$ is obtained as
\begin{equation}\label{eqn:s1_sol}
    S_1(T) = \eta^T S_1(t_0) + \sum_{t=t_0+1}^T \eta^{T-t} N(t) P_1(t,u(t)) \;,
\end{equation}
and correspondingly the total emissions at time $T$, obtained by integrating (\ref{cumE}), is given as
\begin{equation}\label{eqn:cumEsol}
    \mathcal{E}(T) = \mathcal{E}(t_0)  + \sum_{t=t_0+1}^T  M\epsilon_1 \left( \eta^T S_1(t_0) + \sum_{t=t_0+1}^T \eta^{T-t} N(t) P_1(t,u(t)) \right) \;, 
\end{equation}
%
while the total budget for the state, i.e., the sum of annual products of the incentive and the number of EV sales, is
\begin{equation} \label{eqn:IT}
    I(T) = \sum_{t=t_0+1}^T u(t)N_2(t) = \sum_{t=t_0+1}^T u(t)N(t)(1-P_1(t,u(t))) \;.
\end{equation}
Obviously, the higher $u(t)$, the less rapidly $\mathcal{E}(T)$ increases, but at the expense of an increase in $I(T)$.


To find the most appropriate policy roadmap $u(t)$, we formulate an optimal control problem as
\begin{equation} \label{eqn:J}
   \min_{u(t)} I(T) \;,
\end{equation}
subject to state equations (\ref{eqn:S1})-(\ref{cumE}), the terminal condition
\begin{equation} \label{eqn:ebar}
  \mathcal{E}(T) \leq \overline{\mathcal{E}}\;,
\end{equation}
with $\overline{\mathcal{E}}$ as the desired emissions target at horizon $T$, initial condition $\mathcal{E}(t_0) = 0$ and $S_1(t_0)$, and constraints on the control variable as 
\begin{equation}\label{eqn:consU}
    0 \leq u(t) \leq C_v^P(t)\;.
\end{equation}
For the discrete-time system under study, the Hamiltonian is formed as
\begin{equation} \label{eqn:H}
\begin{aligned}
     H(t) =& \; u(t)N(t) (1-P_1(t,u(t))  \; + \\& \lambda(t)  (\eta S_1(t-1) + N(t)(P_1(t,u(t)))) + \\& \nu(t) (\mathcal{E}(t-1) +  M\epsilon_1(\eta S_1(t-1) + N(t)(P_1(t,u(t)))))\;,
\end{aligned}
\end{equation}
so that the optimal policy roadmap, within the control bounds, is found by solving
\begin{equation} \label{eqn:ustar}
\frac{\partial H(t)}{\partial u(t)}  = 0\;.
\end{equation}

The Euler-Lagrange equation yields the dynamics for the adjoint states $\lambda(t)$ and $\nu(t)$,
\begin{equation} \label{eqn:euler}
\begin{aligned}
  \lambda(t-1) &= \frac{\partial H(t)}{\partial S_1(t-1)} = \eta(\lambda(t) +M\epsilon_1 \nu(t))\;,\\
   \nu(t-1) &= \frac{\partial H}{\partial \mathcal{E}(t-1)} = \nu(t)    
\end{aligned}
\end{equation}
with the initial conditions $\lambda(t_0)=\lambda_0$, and $\nu(t_0) = \nu_0$, where $\lambda_0$ and $\nu_0$ are to be found by imposing the terminal constraint (\ref{eqn:ebar}) and the transversality condition $\lambda(T) =0$.
The solution of (\ref{eqn:euler}) is obtained as
\begin{equation}
\begin{aligned}
    \nu(t) &= \nu_{0}\;,\\
      \lambda(t) &= \nu_0 M\epsilon_1 \eta \frac{1-\eta^{T-t+t_0}}{1-\eta}
\end{aligned}
\end{equation}


Assuming the constraints (\ref{eqn:consU}) are satisfied, we obtain the optimal solution $u(t)$ by imposing (\ref{eqn:ustar}), obtaining the implicit equation
\begin{equation} \label{eqn:control}
     \mathcal{Q}(t)\mathcal{R}(t)^{u(t)} / \mathcal{P}(t)+ \ln(\mathcal{R}(t))u(t) = \ln(\mathcal{R}(t))\lambda(t) +   M  \epsilon_1 \ln(\mathcal{R}(t))\nu(t) - 1\;,
\end{equation}
where the sales $N$ or the demand $G$ do not appear explicitly. This transcendental equation can be solved by invoking the Lambert function\footnote{The Lambert function $W(x)$ is the solution of the equation $We^W=x$.}:
\begin{equation}
  u(t) = \frac{a(t)-W\left(-b(t)e^{a(t)}\right)}{c(t)}  \;,
\end{equation}
where $a(t)\triangleq M\epsilon_1 \nu_0 c(t) \left(\frac{\eta^{T-(t-t_0)}-\eta^{-1}}{1-\eta^{-1}}\right) - 1$, $b(t) \triangleq -\mathcal{Q}(t)/\mathcal{P}(t)$, and $c(t) \triangleq \ln(\mathcal{R}(t))$.

By substituting $u$ into (\ref{eqn:cumEsol})--(\ref{eqn:IT}), the target emission level and the total incentive are given by
\begin{equation} \label{eqn:ETlam}
\mathcal{E}(T|\nu_0) =   \sum_{t=t_0+1}^T  M\epsilon_1 \left( \eta^T S_1(t_0) + \sum_{t=t_0+1}^T \eta^{T-t} N(t) \left(1+W(t)\right)^{-1} \right) \;,
\end{equation}
\begin{equation} \label{eqn:ITlam}
  I(T|\nu_0) = \sum_{t=t_0+1}^T N(t) \left(\frac{a(t)-W(t)}{c(t)} \right) W(t) \left(1+W(t)\right)^{-1}  \;,
\end{equation}
where we explicitly denoted the dependence of the solution on the unknown (by now) parameter $\nu_0$, and have written $W(t)$ instead of $W(-b(t)e^{a(t)})$ for the sake of brevity.

It is straightforward to show that
\begin{equation}
    \begin{aligned}
       \frac{\partial }{\partial a} \left( \frac{1}{1+W}\right)    &= -\frac{W} {(1+W)^3}\:, \\
       \frac{\partial }{\partial a} \left( \frac{aW-W^2}{1+W}\right) &= \frac{W(1+a)}{(1+W)^3}\;.
    \end{aligned}
\end{equation}
%
and that, consequently, $\mathcal{E}(T|\nu_0)$ decreases and $I(T|\nu_0)$ increases monotonically with $\nu_0$. Since $I(T)$ is to be minimized, and the constraint (\ref{eqn:ebar}) to be satisfied, it is clear that the optimal solution is obtained for the value of $\nu_0$ such that the equality sign in (\ref{eqn:ebar}) holds, i.e., $\mathcal{E}(T|\nu_0)=\overline{\mathcal{E}}$.
Unfortunately, it is not possible to find this value in closed form. However, it can be found numerically using a shooting algorithm. 
 

Other sensitivities of $\mathcal{E}(T)$, $I(T)$ with respect to key parameters can be evaluated analytically. 
%
In summary the model of this section is useful to give some insight into the nature of the solutions to our backcasting problem, and to show the influence of some factors and parameters. However, it is still too simple to capture essential phenomena, so that its quantitative predictions may be inaccurate.
Therefore we now move on to a more realistic model.

\subsection{Modeling and backcasting: Full-order model}\label{sec:full-order}


When we reintroduce a variable mileage $M(t)$ and $A>1$ age classes, the model becomes more complex, since the survival rate $\eta_{a}$ and the emission factor $\epsilon_{va}$ are now functions of age (too), and so is the old vehicle stock $O_{va}$.

Equation (\ref{eqn:summ}) is replaced by
\begin{equation}
    E(t) = \sum_{a} \epsilon_{1a}(t) M(t) S_{1a}(t)\;,
\end{equation}
equation (\ref{eqn:S1}) by $2 \times(A+1)$ equations
\begin{equation}
\begin{aligned}
S_{v0}(t) &= P_v(t, u(t)) N(t,\textbf{s}(t-1))\\
S_{va}(t) &= \eta_a S_{v,a-1}(t-1), \quad \forall a = 1,\ldots,A-1 \;,\\
S_{vA}(t) &= \eta_A S_{v,A-1}(t-1)+\eta_A S_{vA}(t-1)\;,
\end{aligned}
\end{equation}
and (\ref{cumE}) by 
\begin{equation}
    \mathcal{E}(t) = \mathcal{E}(t-1) + \sum_{a} M(t) \epsilon_{1a}(t) S_{1a}(t) \;,
\end{equation}

where we write $\textbf{s}=[S_{v0},\ldots,S_{vA}]$, $\forall v$, for the stocks state vector. 
The explicit form of the total sales $N(t)$ is
\begin{equation}
 N(t,\textbf{s}(t-1)) = \frac{G(t)}{M(t)} - \sum_{v,a=1}^{A-1} \eta_a S_{v,a-1}(t-1) + \eta_A S_{v,A-1}(t-1) + \eta_A S_{vA}(t-1)\;,
\end{equation}
while $P_1(t,u(t))$ is still given by (\ref{eqn:P}) and (\ref{eqn:p})--(\ref{eqn:R}).


The objective function and the terminal condition for the general OCP are the same as in (\ref{eqn:J})--(\ref{eqn:ebar}), but now the cost function is also a function of the state, due to the age dependency of the survival rates. 

The Hamiltonian is formed as

\begin{equation}
    \begin{aligned}
        H(t) & = u(t) (1 - P_1(u(t))) N(t, \textbf{s}(t-1)) \\
        & + \sum_{v} \lambda_{v0}(t) P_v(u(t)) N(t,\textbf{s}(t-1)) \\
        & + \sum_{v,a=1}^{A-1} \lambda_{va}(t) \eta_a S_{v,a-1}(t-1) \\
        & + \sum_v \lambda_{vA}(t) \eta_A \bigl( S_{v,A-1}(t-1) + S_{vA}(t-1) \bigr) \\
        & + \nu(t) \Bigl( \mathcal{E}(t-1) +  M(t)\epsilon_{10}(t)P_1(u(t))N(t,\textbf{s}(t-1)) \\
        & \qquad \quad  +  \sum_{a=1}^{A-1} M(t)\epsilon_{1a}(t) \eta_a S_{1a-1}(t-1)  \\
        & \qquad  \quad + M(t) \epsilon_{1A}(t)\eta_A(t) \bigl(S_{1A-1}(t-1) + S_{1A}(t-1) \bigr) \Bigr)
    \end{aligned}
\end{equation}
which replaces (\ref{eqn:H}), where $\bm{\lambda}=\{\lambda_{v0},\ldots,\lambda_{vA}\}$, 
$\forall v$ and $\nu(t)$ are the adjoint states.

The first-order optimality conditions (\ref{eqn:ustar})--(\ref{eqn:euler}) now write
\begin{equation}
  \frac{\partial{H(t)}}{\partial{u(t)}} = N(t,\textbf{s}(t-1))\left(1 - P_1 + \frac{\partial{P_1}}{\partial{u(t)}}\left(M(t)\epsilon_{10}(t) \nu(t) + \lambda_{10}(t) - \lambda_{20}(t) - u(t)\right) \right)= 0 \;,
\end{equation}

\begin{equation}
    \begin{aligned}
        \lambda_{va}(t-1) = \frac{\partial{H}(t)}{\partial{S_{va}(t-1)}} &= -\eta_{a+1} \left(u(t)\left(1-P_1\right) + \lambda_{10}(t)P_1 + \lambda_{20}(t)(1-P_1) \right) + \\ &+\lambda_{v,a+1}(t)\eta_{a+1} + M(t) \nu(t) \eta_{a+1} \bigl(\delta_{v1} \epsilon_{1,a+1}(t) - P_1 \epsilon_{10}(t)\bigr) \\&\quad a=0,\ldots,A-1 \;,
    \end{aligned}
\end{equation}

\begin{equation}
\begin{aligned}
      \lambda_{vA}(t-1) = \frac{\partial{H}(t)}{\partial{S_{vA}(t-1)}} =& -\eta_A\left( u(t)(1-P_1) + \lambda_{10}(t) P_1 + \lambda_{20}(t)(1-P_1)\right) + \lambda_{vA}(t) \eta_{A} \\& + M(t) \nu(t) \eta_{A} \bigl(\delta_{v1} \epsilon_{1,A}(t) - P_1 \epsilon_{10}(t)\bigr) \;,
\end{aligned} 
\end{equation}

\begin{equation}
     \nu(t-1) = \frac{\partial H}{\partial \mathcal{E}(t-1)} = \nu(t)\;,
\end{equation}
where we have omitted the dependence of $P_1$ on time and control for the sake of brevity, and $\delta$ represents the Kronecker delta function (that is, 1 if $v=1$ and 0 else). 

This nonlinear system of differential equations cannot be solved analytically as it was the case with (\ref{eqn:ustar})--(\ref{eqn:euler}). Therefore, numerical methods must be used. 
In this work we have used the trust-region constrained algorithm (\textit{trust-constr}) within the Python \textit{scipy.optimize} package.

\section{Results and discussion}\label{sec:results}

\subsection{Case study: Data}




As a case study we consider Metropolitan France to illustrate the proposed backcasting approach. 
For this purpose, the analysis considers a time horizon from $t_0$ = 2022 to $T$ = 2050, with a time step of one year. The CO$_2$ target $\overline{\mathcal{E}}$ is set by forecasting a reference scenario in which a constant incentive (IC) of 5 k\euro~ is provided for each EV purchased, see below.


The parameters of the model are tuned using historical data from the sources listed in Table~\ref{tab:data}.

\begin{table}[t!]
    \centering
    \caption{Data sources}
    \begin{tabular}{|c|c|p{5.5cm}|}
        \hline
        Index & Parameter & Web link \\
        \hline
        1 & $s_{voa}$ & \url{www.statistiques.developpement-durable.gouv.fr/parc-et-circulation-des-vehicules-routiers} \\
        2 & $\epsilon_{1,0}$ & \url{carlabelling.ademe.fr/chiffrescles/r/evolutionTauxCo2} \\
        3 & $e_v$ & \url{www.citepa.org/fr/secten} \\
        4 & $\dot{\chi}$ & \url{www.statistiques.developpement-durable.gouv.fr/immatriculation-des-vehicules-routiers}\\
        \hline
    \end{tabular}
    \label{tab:data}
\end{table}

As for the survival rate, the identification was carried out using data source~1. The latter contains the historical stock of passenger cars $s_{voa}(\tau)$ by technology, ownership type (private, $o = 1$ and professional, $o=2$), and age  until 2022.  Neglecting the dependence of the survival rate on vehicle technology and the movement of second-hand vehicles between ownership types, the survival rate is evaluated as 
\begin{equation}
    \eta_a = \frac{\sum\limits_{vo} s_{voa}(2022)}{\sum\limits_{vo} s_{vo,a-1}(2021)}\;.
\end{equation}
%
%
%

The initial values for the passenger car fleet is set using also data source~1 as
\begin{equation}
    S_{va}(t_0) = s_{v1a}(2022)\;\; \text{and} \;\; S_1(t_0) = \sum\limits_a s_{11a}(2022).
\end{equation}

\begin{figure}[h!]
    \centering
    \begin{subfigure}[t]{0.48\columnwidth}
        \centering
        \pgfplotsset{ylabel shift = -.5em} 
        \input{figParc/na}
    \end{subfigure}
    \begin{subfigure}[t]{0.48\columnwidth}
        \centering
        \pgfplotsset{ylabel shift = -.5em} 
        \input{figParc/EmsissionFactor}
    \end{subfigure}

    \begin{subfigure}[t]{0.48\columnwidth}
        \centering
        \pgfplotsset{ylabel shift = -.5em} 
        \input{figParc/Mileage}
    \end{subfigure}
    \begin{subfigure}[t]{0.48\columnwidth}
        \centering
        \pgfplotsset{ylabel shift = -.5em} 
        \input{figParc/Demand}
    \end{subfigure}
    \caption{Model Inputs and Parameters} \label{fig:param}
\end{figure}
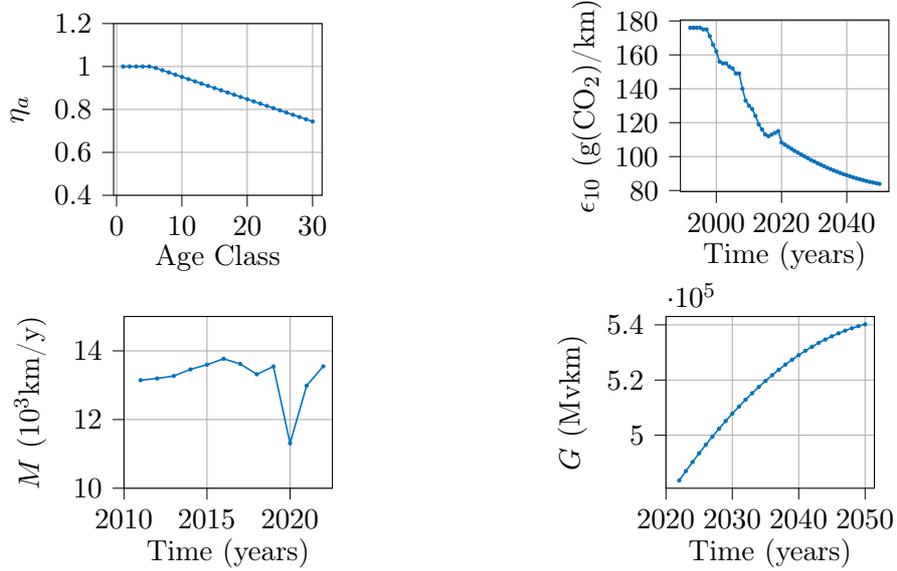

As for the emission factor $\epsilon_{1a}(t)$, which considers only the apparent tailpipe emissions, the identification was carried out using two sources. Data source 2 provides the historical trend ($\tau$ = 1995 to 2020) of average CO$_2$ emissions for newly sold ($a=0$) gasoline and diesel cars. The emission factor for thermal vehicles ($v=1$) for this period is calculated as a weighted average based on the number of newly sold gasoline and diesel vehicles and their respective emission factors. For vehicles sold before 1995, the emission factor is assumed to be at the 1995 level. For the future trend ($\tau$ = 2020 to 2050), (\cite{ITF}) presents the efficiency trajectory of newly sold thermal vehicles in kWh(eq.)/100 km, projecting a 50\% improvement from 2015 to 2050. This trajectory, with the initial value adjusted to be consistent with data source 2, is converted to gCO$_2$/km and approximated by a quadratic function as
\begin{equation}
        \epsilon_{10}(\tau) =  0.01 \cdot (\tau-2020)^2 - 1.27 \cdot (\tau-2020) +  108.2, \quad \tau  \in[2020,2050]\;.
\end{equation}
Given a stock of thermal vehicles by age, their corresponding emission factor $\epsilon_{1a}(t)$ is obtained using the following transformation
\begin{equation}
    \epsilon_{1a} (t) = \epsilon_{10}(t-a)\;.
\end{equation}
The emission factor for EVs is set to zero (i.e., $\epsilon_{2a} = 0$).
The resulting values of the survival rate by vehicle age $\eta_a$ and the emission factor of newly sold thermal vehicles $\epsilon_{10}(\tau)$ are shown in Fig.~\ref{fig:param}.

For the annual mileage $M(t)$, data source 3 provides the annual CO$_2$ emissions $e_1(\tau)$ in France. Assuming that $M(t)$ does not depend on the vehicle type, it is evaluated as 
\begin{equation}\label{eqn:M}
    M(\tau) = \frac{e_1(\tau)}{\sum\limits_{oa} s_{1oa}(\tau)\epsilon_{10}(\tau - a)},\; \tau \in [2011, 2022]\;.
\end{equation}
%
Using (\ref{eqn:M}), Fig.~\ref{fig:param} shows the annual mileage $M(\tau)$, with a dip during the COVID-19 pandemic in 2020. Overall, $M(t)$ is approximated by its average of 13,500~km/y. 

The parameters and determinants related to the logit model in (\ref{eqn:U}) 
are taken from 
(\cite{ADEME}). The transport demand $G(t)$ (vkm), shown in Fig.~\ref{fig:param}, is also taken from that work.

\begin{figure}[h!]
    \centering
    \begin{subfigure}[t]{0.48\columnwidth}
        \centering
        \pgfplotsset{ylabel shift = -.5em} 
        \input{figParc/AdoptionCoeff.tex}
    \end{subfigure}
    \begin{subfigure}[t]{0.48\columnwidth}
        \centering
        \pgfplotsset{ylabel shift = -.5em} 
        \input{figParc/Infra.tex}
    \end{subfigure}
    
    \begin{subfigure}[t]{0.48\columnwidth}
        \centering
        \pgfplotsset{ylabel shift = -.5em} 
        \input{figParc/Costs.tex}
    \end{subfigure}
    \begin{subfigure}[t]{0.48\columnwidth}
        \centering
        \pgfplotsset{ylabel shift = -.5em} 
        \input{figParc/Costs2.tex}
    \end{subfigure}
%
\caption{Determinants used in the logit model.} \label{fig:Exo}
\end{figure}
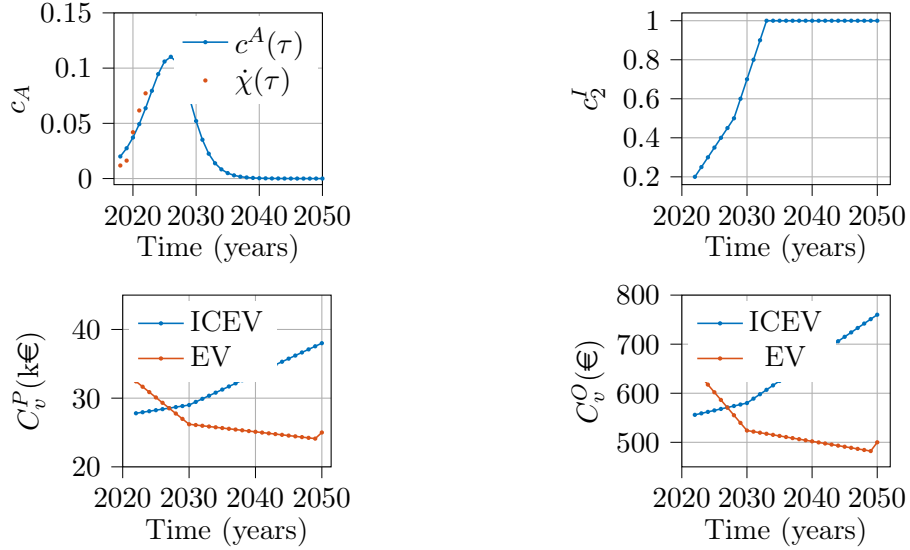

The adoption coefficient $c^A(t)$ is the solution of the Bass model (normalized to market share) 
\begin{equation}\label{eqn:Bass}
\begin{aligned}
    c^A(\tau) = \frac{d}{d\tau}\chi(\tau) = (p + q \chi(\tau))(1-\chi(\tau)) \;,
\end{aligned}
\end{equation}
where 
$p$ and $q$ are the coefficients of innovation and imitation, respectively. The values of $p$ and $q$ are adjusted to match the annual EV sales from 2018 to 2022, as provided by data source 4. Figure~\ref{fig:Exo} and Table~\ref{tab:logit} show the different determinants and parameters used in (\ref{eqn:U}).

\begin{table}[h!]
    \centering
    \caption{List of model parameters.}   
    \label{tab:logit}
    \begin{tabular}{ccc}
        Attributes & ICEV & EV \\ \hline
        Purchase Cost ($p^P$) & -0.3 & -0.3\\
        Operating Cost ($p^O$) & -0.15 & -0.15\\
        Infrastructure Cost ($p^I$) & - & -0.3\\
        $\mu$ & \multicolumn{2}{c}{6.75} \\ 
        $p$ & \multicolumn{2}{c}{0.02} \\ 
        $q$ & \multicolumn{2}{c}{0.4} \\ \hline
    \end{tabular}
    \label{tab:param}
\end{table}

\subsection{Reference Scenarios}

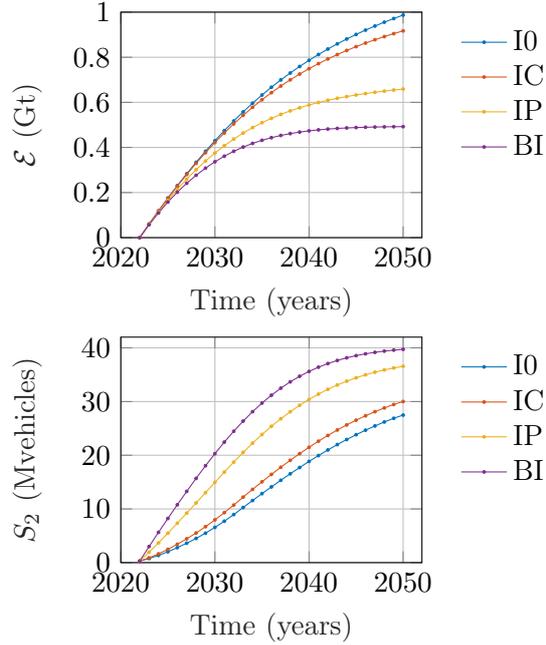
\begin{figure}[h!] 
    \centering
    \begin{subfigure}[t]{\columnwidth}
        \centering
        \input{figParc/caseStudies_cumulEmissions_M2}
    \end{subfigure}
    \hspace{1cm}
    \begin{subfigure}[t]{\columnwidth}
        \centering
        \input{figParc/caseStudies_Stock_cumE_M2}
    \end{subfigure}
    \caption{Reference Scenarios with the full model: CO$_2$ emissions (top) and EV stock (bottom) as a function of time.}\label{fig:refs}
\end{figure}

In addition to the optimal and IC scenarios, we analyze three different policies, namely,
\begin{itemize}
    \item No incentive (I0), $u(t) \equiv 0$ 
    \item Incentive covering the entire EV price (IP), $u(t) = C_2^P(t)$
    \item Ban on ICEV sales from $t_0$ (BI), $N_{2}(t) = N(t)$
\end{itemize}
The corresponding curves of cumulative CO$_2$ emissions and EV stock are shown in Fig.~\ref{fig:refs}. It is clear that the stricter the policy, the faster the increase in EV stock, along with a decrease in cumulative emissions. With the BI scenario, the ICEV stock is virtually depleted by 2050, and consequently the growth of cumulative CO$_2$ emissions ceases by that  year. The $\mathcal{E}(T)$ and $I(T)$ values for the different policy scenarios are given in Table~\ref{tab:res}.

\begin{table}[h!]
    \centering
    \caption{Reference and optimal scenario}    \label{tab:res}
    \begin{tabular}{c|ccccc}
    \hline
              Output  & I0  & IC  & IP  & BI  & Optimal  \\
    \hline
       $\mathcal{E}(T)$ (Gt) &   0.98 & 0.91 & 0.65  & 0.49 & 0.91 \\
         $I(T)$ (G\euro) & 0 & 215 & 1497 & - & 95.3 \\
         \hline
    \end{tabular}
\end{table}

\subsection{Backcasting results}

The IC scenario predicts a total of 0.91~Gt of CO$_2$ in 2050. Accordingly, $\overline{\mathcal{E}}$ is first set to this value to obtain a comparable optimal incentive law. 

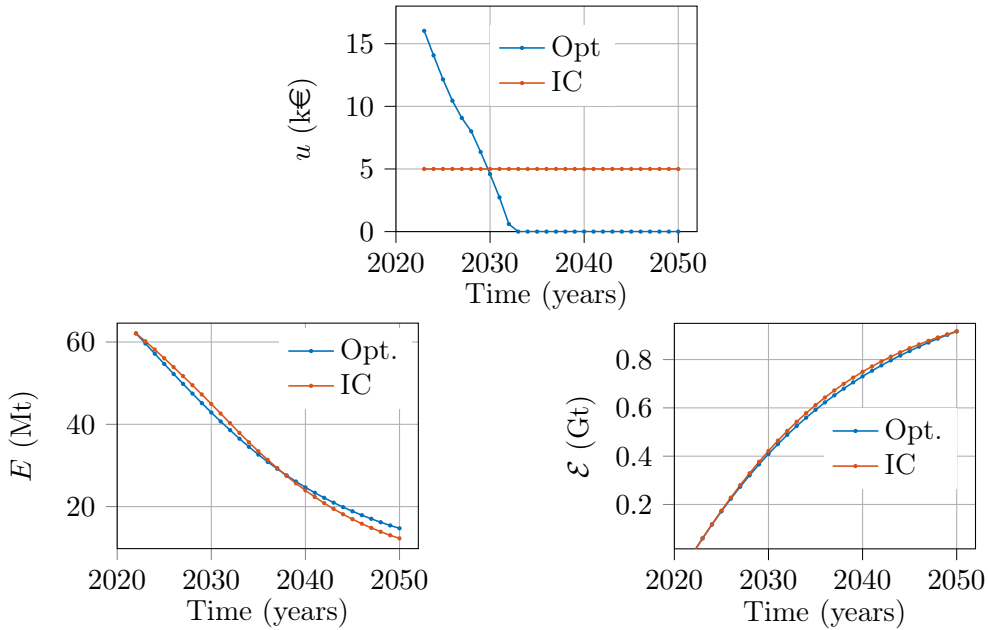
\begin{figure}[h!]
    \centering
    \begin{subfigure}[t]{\columnwidth}
        \centering
         \pgfplotsset{ylabel shift = -.5em} 
        \input{figParc_aac/OptIncentive}
    \end{subfigure}
    
    \begin{subfigure}[t]{.48\columnwidth}
        \centering
         \pgfplotsset{ylabel shift = -.5em} 
        \input{figParc_aac/Emission5KOpt}
    \end{subfigure}
    \begin{subfigure}[t]{.48\columnwidth}
        \centering
        \pgfplotsset{ylabel shift = -.5em} 
        \input{figParc/cumulEmission_M2_heart}
    \end{subfigure}
    \caption{Incentive (top), Yearly Emission (bottom left), and Cumulative Emission (bottom right) Profile} \label{fig:opt}
\end{figure}

The incentive law and emissions curves are shown in Fig.~\ref{fig:opt}.
The incentive law (top panel) exhibits a monotonically decreasing behavior, starting with an incentive close to 50~$\%$ of $C_2^P(t_0)$ until being zero at a certain year. 
Intuitively, such behavior is optimal within the assumptions of the model, as it incentivizes early EV purchases, thereby increasing the EV share in the vehicle stock and reducing the growth of CO$_2$  emissions. 
This effect can be observed in the curves of EV sales and stock, as well as in the emissions curve.

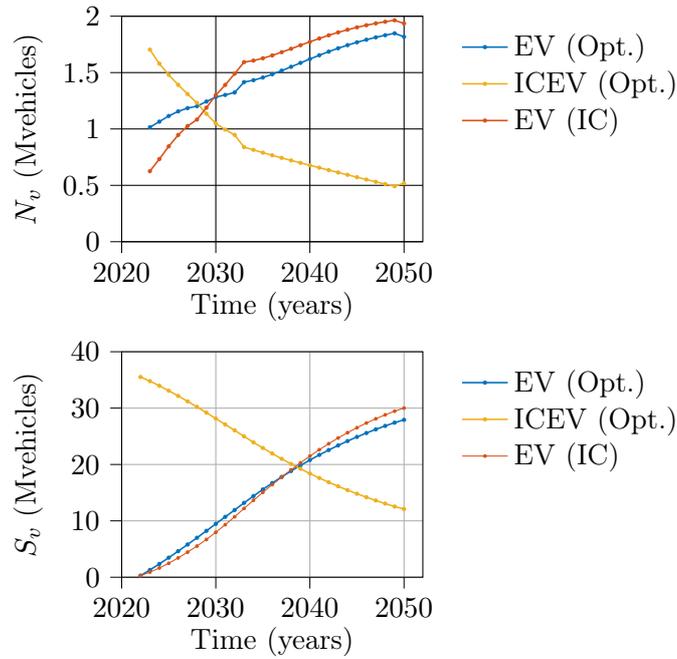
\begin{figure}[h!]
    \centering
    \begin{subfigure}[t]{\columnwidth}
        \centering
        \input{figParc_aac/Sales5KOpt}
    \end{subfigure}
    \hspace{1cm}
    \begin{subfigure}[t]{\columnwidth}
        \centering
        \input{figParc_aac/Stock5KOpt}
    \end{subfigure}
    \caption{Vehicle Sales (top) and Vehicle Stock (bottom) }\label{fig:salesStock}
\end{figure}

The vehicle sales and stock curves for the optimal and IC scenarios are shown in Fig.~\ref{fig:salesStock}.
Compared to the IC law, EV sales are higher in the early years and begin to decrease at 2030.  
Clearly, in both scenarios, the ICEV stock ($v=1$) decreases while the EV stock ($v=2$) increases over time, both exhibiting an S-shaped curve suggesting a variable rate. However, towards the end, the EV stock in the optimal scenario is lower compared to the IC scenario because the optimal strategy focuses on meeting the set final emission constraint. The final cumulative CO$_2$ emissions for both scenarios are the same by construction (constraint at $T$ imposed on the optimal scenario). 

Consequently, Fig.~\ref{fig:opt}, the optimal yearly emissions curve (bottom-left panel) decreases rapidly in the early years and the cumulative emissions curve (bottom-right panel) exhibits a slower growth. 
Overall, we obtain a total expenditure of $I(T)= 95.3$ G\euro, i.e., a reduction of about 55\% compared to the IC reference scenario, see Table~\ref{tab:res}.    
%


\begin{figure}
    \centering
    \input{figParc/Pareto}
    \caption{Pareto Optimal - Incentive vs. Cumulative Emission}
    \label{fig:Pareto}
\end{figure}
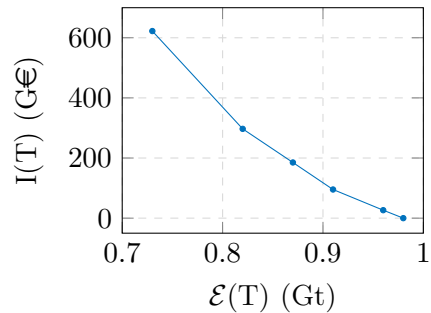

Letting the target on $\mathcal{E}(T)$ vary, we obtain a Pareto optimal trade-off between cumulative CO$_2$ and the total incentive, shown  in Fig.~\ref{fig:Pareto}. The rightmost point indicates the cumulative CO$_2$ emissions at $T$ without any incentive from the state (scenario I0). Conversely, the leftmost point corresponds to the cumulative CO$_2$ emissions at $T$, beyond which the incentives for EVs would exceed their purchase price. It can be seen that the total incentive increases non-linearly with emissions reduction, indicating diminishing returns.
\section{Conclusions}\label{sec:conclusion}

This study represents a first attempt to develop a backcasting methodology to identify the optimal policy roadmaps in transportation systems. 
The analysis focused on a passenger car fleet subsystem, describing its evolution and associated emissions, with the monetary incentive to purchase electric vehicles as the control input. The optimal incentive roadmap was derived by formulating an optimal control problem with the objective of minimizing the government budget while achieving a desired CO$_2$ target. A quantitative case study applied to Metropolitan France was conducted to illustrate the backcasting approach. 

Further research can improve the backcasting paradigm in several ways. Refinements to the fleet model could include regional disaggregation within Metropolitan France, additional vehicle types (e.g., gasoline, diesel, hybrid) and transport modes (e.g., bicycles, rail). Mileage assumptions could be refined by accounting for variation by vehicle type and user profile as in (\cite{ITF}). The survival rate, currently based only on natural obsolescence, could also include factors such as Low Emission Zones (LEZ), which may accelerate vehicle turnover. In addition, the choice of LEZ implementation could be optimized.  Finally, to capture the changes in operating costs and vehicle ownership, the demand for passenger cars, treated here as exogenous, could be replaced by a demand model that predicts vkm by mode and zone.

\section*{Acknowledgements}

This research benefited from state aid managed by the \textit{Agence Nationale de la Recherche (ANR)}, under France 2030, within the project FORBAC bearing the reference ANR-23-PEMO-0002. The authors would like to acknowledge the useful discussions with Dr. Benoit Cheze (IFPEN) and the partners of the FORBAC project. 

\bibliography{parc_aac.bib}

\end{document}

%% file: figParc/na.tex
\definecolor{mycolor1}{rgb}{0.00000,0.44700,0.74100}%
\definecolor{mycolor2}{rgb}{0.85000,0.32500,0.09800}%
\definecolor{mycolor3}{rgb}{0.92900,0.69400,0.12500}%
\definecolor{mycolor4}{rgb}{0.49400,0.18400,0.55600}%
\begin{tikzpicture}

\definecolor{darkgray176}{RGB}{176,176,176}
\definecolor{steelblue31119180}{RGB}{31,119,180}

\begin{axis}[
/pgf/number format/1000 sep={},
scale = 0.4,
tick align=outside,
tick pos=left,
x grid style={darkgray176},
xlabel={Age Class},
xmajorgrids,
xmin=-0.45, xmax=31.45,
xtick style={color=black},
y grid style={darkgray176},
ylabel={$\eta_a$},
ymajorgrids,
ymin=0.4, ymax=1.2,
ytick style={color=black}
]
\addplot [semithick, mycolor1, mark=*, mark size=.5, mark options={solid}]
table {%
1 1
2 1
3 1
4 1
5 1
6 0.99298212
7 0.98260041
8 0.9722187
9 0.96183699
10 0.95145528
11 0.94107357
12 0.93069186
13 0.92031015
14 0.90992844
15 0.89954673
16 0.88916502
17 0.87878331
18 0.8684016
19 0.85801989
20 0.84763818
21 0.83725647
22 0.82687476
23 0.81649305
24 0.80611134
25 0.79572963
26 0.78534792
27 0.77496621
28 0.7645845
29 0.75420279
30 0.74382108
};
\end{axis}

\end{tikzpicture}

%% file: figParc/EmsissionFactor.tex
\definecolor{mycolor1}{rgb}{0.00000,0.44700,0.74100}%
\definecolor{mycolor2}{rgb}{0.85000,0.32500,0.09800}%
\definecolor{mycolor3}{rgb}{0.92900,0.69400,0.12500}%
\definecolor{mycolor4}{rgb}{0.49400,0.18400,0.55600}%
\begin{tikzpicture}

\definecolor{darkgray176}{RGB}{176,176,176}
\definecolor{steelblue31119180}{RGB}{31,119,180}

\begin{axis}[
/pgf/number format/1000 sep={},
scale = 0.4,
tick align=outside,
tick pos=left,
x grid style={darkgray176},
xlabel={Time (years)},
xmajorgrids,
xmin=1989, xmax=2052.9,
xtick style={color=black},
y grid style={darkgray176},
ylabel={$\epsilon_{10}$ (g(CO$_2$)/km)},
ymajorgrids,
ymin=79.295, ymax=180.605,
ytick style={color=black}
]
\addplot [semithick, mycolor1, mark=*, mark size=.5, mark options={solid}]
table {%
1992 176
1993 176
1994 176
1995 176
1996 175
1997 175
1998 171
1999 166
2000 162
2001 156
2002 155
2003 155
2004 153
2005 152
2006 149
2007 149
2008 140
2009 133
2010 130
2011 128
2012 124
2013 119
2014 116
2015 113
2016 112
2017 113
2018 114
2019 115
2020 108.3
2021 107
2022 105.8
2023 104.6
2024 103.4
2025 102.3
2026 101.2
2027 100.1
2028 99.1
2029 98
2030 97.1
2031 96.1
2032 95.2
2033 94.3
2034 93.4
2035 92.6
2036 91.8
2037 91.1
2038 90.3
2039 89.6
2040 89
2041 88.3
2042 87.7
2043 87.1
2044 86.6
2045 86.1
2046 85.6
2047 85.1
2048 84.7
2049 84.3
2050 83.9
};
\end{axis}

\end{tikzpicture}

%% file: figParc/Mileage.tex
\definecolor{mycolor1}{rgb}{0.00000,0.44700,0.74100}%
\definecolor{mycolor2}{rgb}{0.85000,0.32500,0.09800}%
\definecolor{mycolor3}{rgb}{0.92900,0.69400,0.12500}%
\definecolor{mycolor4}{rgb}{0.49400,0.18400,0.55600}%

\begin{tikzpicture}

\definecolor{darkgray176}{RGB}{176,176,176}
\definecolor{steelblue31119180}{RGB}{31,119,180}

\begin{axis}[
/pgf/number format/1000 sep={},
scale=0.4,
tick align=outside,
tick pos=left,
x grid style={darkgray176},
xlabel={Time (years)},
xmajorgrids,
xmin=2010, xmax=2022.55,
xtick style={color=black},
y grid style={darkgray176},
ylabel={$M$ (10$^3$km/y)},
ymajorgrids,
ymin=10, ymax=15,
ytick style={color=black}
]
\addplot [semithick, mycolor1,  mark=*, mark size=.5, mark options={solid}]
table {%
2011 13.146
2012 13.196
2013 13.269
2014 13.460
2015 13.597
2016 13.769
2017 13.621
2018 13.318
2019 13.544
2020 11.309
2021 12.988
2022 13.550
};
\end{axis}

\end{tikzpicture}

%% file: figParc/Demand.tex
\definecolor{mycolor1}{rgb}{0.00000,0.44700,0.74100}%
\definecolor{mycolor2}{rgb}{0.85000,0.32500,0.09800}%
\definecolor{mycolor3}{rgb}{0.92900,0.69400,0.12500}%
\definecolor{mycolor4}{rgb}{0.49400,0.18400,0.55600}%

\begin{tikzpicture}

\definecolor{darkgray176}{RGB}{176,176,176}
\definecolor{steelblue31119180}{RGB}{31,119,180}

\begin{axis}[
/pgf/number format/1000 sep={},
scale = 0.4,
tick align=outside,
tick pos=left,
x grid style={darkgray176},
xlabel={Time (years)},
xmajorgrids,
xmin=2020, xmax=2051.4,
xtick style={color=black},
y grid style={darkgray176},
ylabel={$G$ (Mvkm)},
ymajorgrids,
ymin=480791.457693442, ymax=543015.799327685,
ytick style={color=black}
]
\addplot [semithick, mycolor1, mark=*, mark size=.5, mark options={solid}]
table {%
2022 483619.836858635
2023 486997.48528258
2024 490274.587218235
2025 493451.142665625
2026 496527.151624724
2027 499502.614095582
2028 502377.530078175
2029 505151.899572478
2030 507825.722578466
2031 510398.999096238
2032 512871.72912572
2033 515243.912666961
2034 517515.549719911
2035 519686.640284596
2036 521757.184361016
2037 523727.18194917
2038 525596.633049058
2039 527365.537660657
2040 529033.895784014
2041 530601.707419081
2042 532068.972565809
2043 533435.691224346
2044 534701.863394616
2045 535867.489076597
2046 536932.568270312
2047 537897.100975762
2048 538761.087192946
2049 539524.526921864
2050 540187.420162492
};
\end{axis}

\end{tikzpicture}

%% file: figParc/AdoptionCoeff.tex
\definecolor{mycolor1}{rgb}{0.00000,0.44700,0.74100}%
\definecolor{mycolor2}{rgb}{0.85000,0.32500,0.09800}%
\definecolor{mycolor3}{rgb}{0.92900,0.69400,0.12500}%
\definecolor{mycolor4}{rgb}{0.49400,0.18400,0.55600}%

\begin{tikzpicture}

\definecolor{darkgray176}{RGB}{176,176,176}
\definecolor{lightgray204}{RGB}{204,204,204}
\definecolor{steelblue31119180}{RGB}{31,119,180}

\begin{axis}[
/pgf/number format/1000 sep={},
scale = 0.4,
legend cell align={left},
legend style={fill opacity=0, draw opacity=1, text opacity=1, draw=none},
tick align=outside,
tick pos=left,
x grid style={darkgray176},
xlabel={Time (years)},
xmajorgrids,
xmin=2017, xmax=2050,
xtick style={color=black},
y grid style={darkgray176},
ylabel={ $c_A$},
ymajorgrids,
ymin=-0.0055125, ymax=0.15,
ytick style={color=black},
yticklabel style={
        /pgf/number format/fixed,
        /pgf/number format/precision=3
},
scaled y ticks=false
]
\addplot [semithick, mycolor1, mark=*, mark size=.5, mark options={solid}]
table {%
2018 0.02
2019 0.02744
2020 0.03713
2021 0.04927
2022 0.06369
2023 0.07946
2024 0.09457
2025 0.10597
2026 0.11025
2027 0.10516
2028 0.09125
2029 0.07201
2030 0.05212
2031 0.03514
2032 0.02247
2033 0.01384
2034 0.00833
2035 0.00494
2036 0.0029
2037 0.0017
2038 0.00099
2039 0.00057
2040 0.00033
2041 0.00019
2042 0.00011
2043 7e-05
2044 4e-05
2045 2e-05
2046 1e-05
2047 1e-05
2048 0
2049 0
2050 0
};
\addlegendentry{$c^A(\tau)$}
\addplot [semithick, mycolor2, mark=*, mark size=.5, mark options={solid}, only marks]
table {%
2018 0.0118069936061286
2019 0.0162375472072179
2020 0.041949571268378
2021 0.0616447665732201
2022 0.0771788740705876
};
\addlegendentry{$\dot{\chi}(\tau)$}
\end{axis}

\end{tikzpicture}

%% file: figParc/Infra.tex
\definecolor{mycolor1}{rgb}{0.00000,0.44700,0.74100}%
\definecolor{mycolor2}{rgb}{0.85000,0.32500,0.09800}%
\definecolor{mycolor3}{rgb}{0.92900,0.69400,0.12500}%
\definecolor{mycolor4}{rgb}{0.49400,0.18400,0.55600}%

\begin{tikzpicture}

\definecolor{darkgray176}{RGB}{176,176,176}
\definecolor{steelblue31119180}{RGB}{31,119,180}

\begin{axis}[
scale=0.4,
/pgf/number format/1000 sep={},
tick align=outside,
tick pos=left,
x grid style={darkgray176},
xlabel={Time (years)},
xmajorgrids,
xmin=2020, xmax=2052,
xtick style={color=black},
y grid style={darkgray176},
ylabel={$c_2^I$},
ymajorgrids,
ymin=0.16, ymax=1.04,
ytick style={color=black}
]
\addplot [semithick, mycolor1, mark=*, mark size=.5, mark options={solid}]
table {%
2022 0.2
2023 0.25
2024 0.3
2025 0.35
2026 0.4
2027 0.45
2028 0.5
2029 0.6
2030 0.7
2031 0.8
2032 0.9
2033 1
2034 1
2035 1
2036 1
2037 1
2038 1
2039 1
2040 1
2041 1
2042 1
2043 1
2044 1
2045 1
2046 1
2047 1
2048 1
2049 1
2050 1
};
\end{axis}

\end{tikzpicture}

%% file: figParc/Costs.tex
\definecolor{mycolor1}{rgb}{0.00000,0.44700,0.74100}%
\definecolor{mycolor2}{rgb}{0.85000,0.32500,0.09800}%
\definecolor{mycolor3}{rgb}{0.92900,0.69400,0.12500}%
\definecolor{mycolor4}{rgb}{0.49400,0.18400,0.55600}%

\begin{tikzpicture}

\definecolor{darkgray176}{RGB}{176,176,176}
\definecolor{darkorange25512714}{RGB}{255,127,14}
\definecolor{steelblue31119180}{RGB}{31,119,180}

\begin{axis}[
/pgf/number format/1000 sep={},
scale = 0.4,
tick align=outside,
tick pos=left,
x grid style={darkgray176},
xlabel={},
xmajorgrids,
xmin=2020, xmax=2051.4,
xtick style={color=black},
xlabel={Time (years)},
y grid style={darkgray176},
ylabel=
{$C_v^P$(k\euro)},
ymajorgrids,
ymin=20.000, ymax=45.000,
ytick style={color=black},
legend pos=north west, 
legend cell align={left},
legend style={fill opacity=0, draw opacity=1, text opacity=1, draw=none},
yticklabel style={
        /pgf/number format/fixed,
        /pgf/number format/precision=4
},
scaled y ticks=false
]
\addplot [semithick, mycolor1, mark=*, mark size=.5, mark options={solid}]
table {%
2022 27.800
2023 27.950
2024 28.100
2025 28.250
2026 28.400
2027 28.550
2028 28.700
2029 28.850
2030 29.000
2031 29.450
2032 29.900
2033 30.350
2034 30.800
2035 31.250
2036 31.700
2037 32.150
2038 32.600
2039 33.050
2040 33.500
2041 33.950
2042 34.400
2043 34.850
2044 35.300
2045 35.750
2046 36.200
2047 36.650
2048 37.100
2049 37.550
2050 38.000
};
\addlegendentry{ICEV}
\addplot [semithick, mycolor2, mark=*, mark size=.5, mark options={solid}]
table {%
2022 32.440
2023 31.660
2024 30.880
2025 30.100
2026 29.320
2027 28.540
2028 27.760
2029 26.980
2030 26.200
2031 26.090
2032 25.980
2033 25.870
2034 25.760
2035 25.650
2036 25.540
2037 25.430
2038 25.320
2039 25.210
2040 25.100
2041 24.990
2042 24.880
2043 24.770
2044 24.660
2045 24.550
2046 24.440
2047 24.330
2048 24.220
2049 24.110
2050 25.000
};
\addlegendentry{EV}
\end{axis}

\end{tikzpicture}

%% file: figParc/Costs2.tex
\definecolor{mycolor1}{rgb}{0.00000,0.44700,0.74100}%
\definecolor{mycolor2}{rgb}{0.85000,0.32500,0.09800}%
\definecolor{mycolor3}{rgb}{0.92900,0.69400,0.12500}%
\definecolor{mycolor4}{rgb}{0.49400,0.18400,0.55600}%

\begin{tikzpicture}

\definecolor{darkgray176}{RGB}{176,176,176}
\definecolor{darkorange25512714}{RGB}{255,127,14}
\definecolor{steelblue31119180}{RGB}{31,119,180}


\begin{axis}[
/pgf/number format/1000 sep={},
scale = 0.4,
tick pos=left,
tick align=outside,
x grid style={darkgray176},
xmin=2020, xmax=2051.4,
xtick pos=left,
xlabel={Time (years)},
xmajorgrids,
xmin=2020, xmax=2052,
xtick style={color=black},
y grid style={darkgray176},
ymajorgrids,
ylabel=
{$C_v^O$(\euro)},
ymin=450, ymax=800,
ytick style={color=black},
legend pos=north west, 
legend style={fill opacity=0, draw opacity=1, text opacity=1, draw=none},
]
\addplot [semithick, mycolor1, mark=*, mark size=.5, mark options={solid}]
table {%
2022 556
2023 559
2024 562
2025 565
2026 568
2027 571
2028 574
2029 577
2030 580
2031 589
2032 598
2033 607
2034 616
2035 625
2036 634
2037 643
2038 652
2039 661
2040 670
2041 679
2042 688
2043 697
2044 706
2045 715
2046 724
2047 733
2048 742
2049 751
2050 760
};
\addlegendentry{ICEV}
\addplot [semithick, mycolor2, mark=*, mark size=.5, mark options={solid}]
table {%
2022 648.8
2023 633.2
2024 617.6
2025 602
2026 586.4
2027 570.8
2028 555.2
2029 539.6
2030 524
2031 521.8
2032 519.6
2033 517.4
2034 515.2
2035 513
2036 510.8
2037 508.6
2038 506.4
2039 504.2
2040 502
2041 499.8
2042 497.6
2043 495.4
2044 493.2
2045 491
2046 488.8
2047 486.6
2048 484.4
2049 482.2
2050 500
};
\addlegendentry{EV}
\end{axis}

\end{tikzpicture}

%% file: figParc/caseStudies_cumulEmissions_M2.tex
%
%
\definecolor{mycolor1}{rgb}{0.00000,0.44700,0.74100}%
\definecolor{mycolor2}{rgb}{0.85000,0.32500,0.09800}%
\definecolor{mycolor3}{rgb}{0.92900,0.69400,0.12500}%
\definecolor{mycolor4}{rgb}{0.49400,0.18400,0.55600}%
\begin{tikzpicture}

\begin{axis}[%
/pgf/number format/1000 sep={},
width=4.521in,
height=3.566in,
at={(0.758in,0.481in)},
scale =.4,
xmin=2020,
xmax=2052,
xlabel style={font=\color{white!15!black}},
xlabel={Time (years)},
ymin=0,
ymax=1,
ylabel style={font=\color{white!15!black}},
ylabel={$\mathcal{E}$ (Gt)},
axis background/.style={fill=white},
xmajorgrids,
ymajorgrids,
legend style={at={(1.1,.97)}, anchor=north west, legend cell align=left, align=left, draw=none, fill opacity=0.8}
]
\addplot [color=mycolor1, mark=*, mark options={solid, mycolor1}, mark size=.5]
  table[row sep=crcr]{
2022    0          \\
2023    0.0603772  \\
2024    0.118943   \\
2025    0.175636   \\
2026    0.230406   \\
2027    0.283235   \\
2028    0.33412    \\
2029    0.383009   \\
2030    0.429849   \\
2031    0.47461    \\
2032    0.517251   \\
2033    0.55772    \\
2034    0.596095   \\
2035    0.632449   \\
2036    0.666847   \\
2037    0.69938    \\
2038    0.730085   \\
2039    0.759054   \\
2040    0.786372   \\
2041    0.812128   \\
2042    0.836393   \\
2043    0.859251   \\
2044    0.880789   \\
2045    0.901092   \\
2046    0.92024    \\
2047    0.938319   \\
2048    0.955396   \\
2049    0.971547   \\
2050    0.986869   \\
};
\addlegendentry{I0}

\addplot [color=mycolor2, mark=*, mark options={solid, mycolor2}, mark size=.5]
  table{%
2022 0.000000000000
2023 0.060183606657653
2024 0.118342953300956
2025 0.174400438543666
2026 0.228293224328208
2027 0.279996184846205
2028 0.329498902109346
2029 0.376745097389054
2030 0.42168360420384
2031 0.464295754804687
2032 0.504562276425067
2033 0.542460585687819
2034 0.578087622412111
2035 0.611533783857183
2036 0.642888565734561
2037 0.672260432438132
2038 0.69971133868415
2039 0.725352958737131
2040 0.749294651807535
2041 0.771642422594029
2042 0.792487378729662
2043 0.811931101134623
2044 0.830075496816103
2045 0.84701934303609
2046 0.862856599324147
2047 0.877681519291866
2048 0.891573266864887
2049 0.904611702687249
2050 0.916893927455687
};
\addlegendentry{IC}

\addplot [color=mycolor3, mark=*, mark options={solid, mycolor3}, mark size=.5]
   table{%
2022    0.0
2023    0.05872236
2024    0.11402357
2025    0.16589833
2026    0.21436485
2027    0.25947774
2028    0.30130238
2029    0.3398913
2030    0.37532569
2031    0.407716
2032    0.43719024
2033    0.46388383
2034    0.48798399
2035    0.50967519
2036    0.52914133
2037    0.54657715
2038    0.56213142
2039    0.57599523
2040    0.58834859
2041    0.59935389
2042    0.60915373
2043    0.61788691
2044    0.62568762
2045    0.63267385
2046    0.63895482
2047    0.64462906
2048    0.64977772
2049    0.65447406
2050    0.65877421
};
\addlegendentry{IP}

 \addplot [color=mycolor4, mark=*, mark options={solid, mycolor4}, mark size=.5]
   table{%
2022    0.0
2023    0.05724569
2024    0.10981006
2025    0.15786314
2026    0.20156905
2027    0.24110258
2028    0.27663082
2029    0.30835014
2030    0.3364795
2031    0.36125127
2032    0.38291123
2033    0.40170626
2034    0.41788711
2035    0.43170164
2036    0.44339445
2037    0.45321385
2038    0.46135799
2039    0.46805951
2040    0.47353156
2041    0.47795919
2042    0.48150151
2043    0.48430424
2044    0.48650353
2045    0.48821148
2046    0.48952832
2047    0.49053789
2048    0.49130558
2049    0.49188533
2050    0.49230446
};
\addlegendentry{BI}

\end{axis}

\end{tikzpicture}%

%% file: figParc/caseStudies_Stock_cumE_M2.tex
%
%
\definecolor{mycolor1}{rgb}{0.00000,0.44700,0.74100}%
\definecolor{mycolor2}{rgb}{0.85000,0.32500,0.09800}%
\definecolor{mycolor3}{rgb}{0.92900,0.69400,0.12500}%
\definecolor{mycolor4}{rgb}{0.49400,0.18400,0.55600}%
\begin{tikzpicture}

\begin{axis}[%
/pgf/number format/1000 sep={},
width=4.521in,
height=3.566in,
at={(0.758in,0.481in)},
scale =.4,
xmin=2020,
xmax=2052,
xlabel style={font=\color{white!15!black}},
xlabel={Time (years)},
ymin=0,
ymax=42,
ylabel style={font=\color{white!15!black}},
ylabel={$S_2$ (Mvehicles)},
axis background/.style={fill=white},
xmajorgrids,
ymajorgrids,
legend style={at={(1.1,0.97)}, 
anchor=north west, 
legend cell align=left, 
align=left, 
draw=none, 
fill opacity=0.8}
]
\addplot [color=mycolor1, mark=*, mark options={solid, mycolor1},, mark size=.5]
  table{%
2022    0.276883
2023    0.76409072
2024    1.34310864
2025    2.02060232
2026    2.7874171
2027    3.62453088
2028    4.51102615
2029    5.48964332
2030    6.56396393
2031    7.71554824
2032    8.95042866
2033    10.27005694
2034    11.56907137
2035    12.84801638
2036    14.10544952
2037    15.33848709
2038    16.54315827
2039    17.71461936
2040    18.84811574
2041    19.93941685
2042    20.98497663
2043    21.98199123
2044    22.92847897
2045    23.82321502
2046    24.66570774
2047    25.45616052
2048    26.19543404
2049    26.88499518
2050    27.48064894
};
\addlegendentry{I0}

\addplot [color=mycolor2, mark=*, mark options={solid, mycolor2}, mark size=.5]
  table[row sep=crcr]{%
2022	0.276883\\
2023	0.90198211267876\\
2024	1.63387948077211\\
2025	2.47794174631282\\
2026	3.42186901063984\\
2027	4.44391861113421\\
2028	5.52124728177362\\
2029	6.69710325793011\\
2030	7.97000666968697\\
2031	9.31230240582858\\
2032	10.7239404782539\\
2033	12.2002755811868\\
2034	13.6441437964197\\
2035	15.0546221112408\\
2036	16.4291370157325\\
2037	17.7640184007189\\
2038	19.0548684113387\\
2039	20.2967919404005\\
2040	21.4852832249128\\
2041	22.6165959604329\\
2042	23.6878800947605\\
2043	24.6972043125328\\
2044	25.6436017441817\\
2045	26.5269593070818\\
2046	27.3479553979408\\
2047	28.1079858433601\\
2048	28.8090953418799\\
2049	29.4539002713503\\
2050	30.0061194771407\\
};
\addlegendentry{IC}

\addplot [color=mycolor3, mark=*, mark options={solid, mycolor3}, mark size=.5]
  table{%
2022    0.276883
2023    1.9427586
2024    3.6792428
2025    5.48524921
2026    7.34119965
2027    9.22532531
2028    11.11672762
2029    13.03315974
2030    14.95449197
2031    16.85461972
2032    18.7185467
2033    20.53205211
2034    22.25058287
2035    23.86984267
2036    25.38640476
2037    26.79790734
2038    28.10314811
2039    29.3020929
2040    30.39611333
2041    31.38788615
2042    32.28128786
2043    33.08116922
2044    33.79318506
2045    34.42347718
2046    34.97845884
2047    35.46463138
2048    35.88845265
2049    36.25623718
2050    36.56490913
};
\addlegendentry{IP}

\addplot [color=mycolor4, mark=*, mark options={solid, mycolor4}, mark size=.5]
  table{%
2022    0.276883
2023    2.99451343
2024    5.63727769
2025    8.23111527
2026    10.77567125
2027    13.26628132
2028    15.6941916
2029    18.04433181
2030    20.29988585
2031    22.44578992
2032    24.46930826
2033    26.36030838
2034    28.11145226
2035    29.718253
2036    31.17898297
2037    32.49455002
2038    33.66828737
2039    34.70566353
2040    35.6139858
2041    36.40192242
2042    37.07917976
2043    37.65607814
2044    38.14322271
2045    38.5510591
2046    38.88958143
2047    39.16811683
2048    39.39520371
2049    39.57853616
2050    39.72488651

};
\addlegendentry{BI}

\end{axis}

\end{tikzpicture}%

%% file: figParc_aac/OptIncentive.tex
\definecolor{mycolor1}{rgb}{0.00000,0.44700,0.74100}%
\definecolor{mycolor2}{rgb}{0.85000,0.32500,0.09800}%
\definecolor{mycolor3}{rgb}{0.92900,0.69400,0.12500}%
\definecolor{mycolor4}{rgb}{0.49400,0.18400,0.55600}%

\begin{tikzpicture}

\definecolor{darkgray176}{RGB}{176,176,176}
\definecolor{lightgray204}{RGB}{204,204,204}
\definecolor{steelblue31119180}{RGB}{31,119,180}

\begin{axis}[
/pgf/number format/1000 sep={},
width=4.521in,
height=3.566in,
at={(0.758in,0.481in)},
scale=0.4,
legend cell align={left},
legend style={
  fill opacity=0,
  draw opacity=1,
  text opacity=1,
  at= {(axis cs: 2030,10)},
  anchor=south west,
  draw=none
},
tick align=outside,
tick pos=left,
x grid style={darkgray176},
xlabel={Time (years)},
xmajorgrids,
xmin=2020, xmax=2052,
xtick style={color=black},
y grid style={darkgray176},
ylabel={$u$ (k\euro) },
ymajorgrids,
ymin=0, ymax=18,
ytick style={color=black}
]
\addplot [semithick, mycolor1, mark=*, mark size=.5, mark options={solid}]
table {%

2023 16.027071352215
2024 14.06938781574
2025 12.15172890585
2026 10.44599065068
2027 9.072187285755
2028 8.0056819002
2029 6.354953454225
2030 4.5890303616
2031 2.732734600561
2032 0.604945902374
2033 5.89397748981e-09
2034 1.143075510108e-08
2035 2.031528948005e-10
2036 6.3794688345e-09
2037 5.25227425541e-09
2038 8.80376114192e-11
2039 3.65187944961e-09
2040 1.70971813564e-09
2041 2.858319209011e-09
2042 2.092646033556e-09
2043 5.84913549451e-10
2044 2.146455329164e-09
2045 2.08022528754e-09
2046 1.692750388456e-10
2047 1.276257252809e-09
2048 7.0476540015e-10
2049 2.293180393302e-10
2050 9.859974642e-11

};
\addlegendentry{Opt}

\addplot [semithick, mycolor2, mark=*, mark size=.5, mark options={solid}]
table {%
2023 5
2024 5
2025 5
2026 5
2027 5
2028 5
2029 5
2030 5
2031 5
2032 5
2033 5
2034 5
2035 5
2036 5
2037 5
2038 5
2039 5
2040 5
2041 5
2042 5
2043 5
2044 5
2045 5
2046 5
2047 5
2048 5
2049 5
2050 5
};
\addlegendentry{IC}

\end{axis}

\end{tikzpicture}

%% file: figParc_aac/Emission5KOpt.tex
\definecolor{mycolor1}{rgb}{0.00000,0.44700,0.74100}%
\definecolor{mycolor2}{rgb}{0.85000,0.32500,0.09800}%
\definecolor{mycolor3}{rgb}{0.92900,0.69400,0.12500}%
\definecolor{mycolor4}{rgb}{0.49400,0.18400,0.55600}%

\begin{tikzpicture}

\definecolor{darkgray176}{RGB}{176,176,176}
\definecolor{lightgray204}{RGB}{204,204,204}
\definecolor{steelblue31119180}{RGB}{31,119,180}

\begin{axis}[
/pgf/number format/1000 sep={},
width=4.521in,
height=3.566in,
at={(0.758in,0.481in)},
scale=0.4,
legend cell align={left},
legend style={fill opacity=0, draw opacity=1, text opacity=1, draw=none},
tick align=outside,
tick pos=left,
x grid style={darkgray176},
xlabel={Time (years)},
xmajorgrids,
xmin=2020, xmax=2052,
y grid style={darkgray176},
ylabel={$E$ (Mt)},
ymajorgrids,
ymin=9.79146907119933, ymax=64.5880943730381,
]
\addplot [semithick, mycolor1, mark=*, mark options={dashed, mycolor1}, mark size=.5, mark options={solid}]
table {%
2022 62.0973386775
2023 59.6369273092608
2024 57.1515060019351
2025 54.6786129429341
2026 52.2284823100033
2027 49.8241624410995
2028 47.4653029682664
2029 45.1409349468949
2030 42.8674930762398
2031 40.6816046421603
2032 38.5899245895381
2033 36.5040089862257
2034 34.5256276910336
2035 32.6470704667808
2036 30.865066296883
2037 29.1971991201004
2038 27.5899303856052
2039 26.0877697603396
2040 24.6833631878252
2041 23.3707392667316
2042 22.1288935954514
2043 20.9668524921619
2044 19.8804295251482
2045 18.8665690866638
2046 17.9177302468571
2047 17.0357915775675
2048 16.2045649240206
2049 15.4275600754105
2050 14.7298093987136
};
\addlegendentry{Opt.}

\addplot [semithick,color=mycolor2, mark=*, mark options={dashed, mycolor2}, mark size=.5]
  table[row sep=crcr]{%
2022	62.0973386775\\
2023	60.1836066576528\\
2024	58.1593466433033\\
2025	56.0574852427096\\
2026	53.8927857845427\\
2027	51.7029605179969\\
2028	49.5027172631407\\
2029	47.2461952797078\\
2030	44.9385068147861\\
2031	42.6121506008471\\
2032	40.2665216203798\\
2033	37.8983092627524\\
2034	35.6270367242917\\
2035	33.446161445072\\
2036	31.3547818773782\\
2037	29.3718667035712\\
2038	27.4509062460178\\
2039	25.6416200529808\\
2040	23.9416930704035\\
2041	22.347770786495\\
2042	20.8449561356326\\
2043	19.4437224049606\\
2044	18.1443956814805\\
2045	16.9438462199866\\
2046	15.8372562880573\\
2047	14.8249199677192\\
2048	13.8917475730213\\
2049	13.0384358223612\\
2050	12.2822247684378\\
};
\addlegendentry{IC}


\end{axis}

\end{tikzpicture}

%% file: figParc/cumulEmission_M2_heart.tex
\begin{tikzpicture}

\definecolor{mycolor1}{rgb}{0.00000,0.44700,0.74100}%
\definecolor{mycolor2}{rgb}{0.85000,0.32500,0.09800}%
\definecolor{mycolor3}{rgb}{0.92900,0.69400,0.12500}%
\definecolor{mycolor4}{rgb}{0.49400,0.18400,0.55600}%
\definecolor{darkgray176}{RGB}{176,176,176}
\definecolor{lightgray204}{RGB}{204,204,204}
\begin{axis}[
/pgf/number format/1000 sep={},
width=4.521in,
height=3.566in,
at={(0.758in,0.481in)},
scale=0.4,
legend cell align={left},
legend style={fill opacity=0, draw opacity=1, text opacity=1, draw=none, at= {(axis cs: 2050,.6)},},
tick align=outside,
tick pos=left,
x grid style={darkgray176},
xlabel={Time (years)},
xmajorgrids,
xmin=2020, xmax=2052,
xtick style={color=black},
y grid style={darkgray176},
xtick pos=left,
xtick style={color=black},
y grid style={darkgray176},
ylabel={$ \mathcal{E}$ (Gt)},
ymajorgrids,
ymin=0.0162526423047157, ymax=0.95,
ytick pos=left,
]
\addplot [semithick, mycolor1, mark=*, mark size=0.5, mark options={solid}]
table {%
2022 0.000000000000
2023 0.059636927309261
2024 0.116788433311196
2025 0.17146704625413
2026 0.223695528564133
2027 0.273519690905233
2028 0.320985993973499
2029 0.366126928920394
2030 0.408994421996634
2031 0.449676026638794
2032 0.488265951228332
2033 0.524769960214558
2034 0.559295587905592
2035 0.591942658372372
2036 0.622807724669255
2037 0.652004923789356
2038 0.679594854174961
2039 0.705682623935301
2040 0.730366987123126
2041 0.753737726389858
2042 0.775866619985309
2043 0.796833472477471
2044 0.816713901002619
2045 0.835580470089283
2046 0.85349820033614
2047 0.870534991913708
2048 0.886739556837728
2049 0.902167116913139
2050 0.916896926311852

};
\addlegendentry{Opt.}
\addplot [semithick, mycolor2,  mark=*, mark size=0.5, mark options={solid}]
table {%
2022 0.000000000000
2023 0.060183606657653
2024 0.118342953300956
2025 0.174400438543666
2026 0.228293224328208
2027 0.279996184846205
2028 0.329498902109346
2029 0.376745097389054
2030 0.42168360420384
2031 0.464295754804687
2032 0.504562276425067
2033 0.542460585687819
2034 0.578087622412111
2035 0.611533783857183
2036 0.642888565734561
2037 0.672260432438132
2038 0.69971133868415
2039 0.725352958737131
2040 0.749294651807535
2041 0.771642422594029
2042 0.792487378729662
2043 0.811931101134623
2044 0.830075496816103
2045 0.84701934303609
2046 0.862856599324147
2047 0.877681519291866
2048 0.891573266864887
2049 0.904611702687249
2050 0.916893927455687
};
\addlegendentry{IC}
\end{axis}

\end{tikzpicture}

%% file: figParc_aac/Sales5KOpt.tex
\begin{tikzpicture}

\definecolor{mycolor1}{rgb}{0.00000,0.44700,0.74100}%
\definecolor{mycolor2}{rgb}{0.85000,0.32500,0.09800}%
\definecolor{mycolor3}{rgb}{0.92900,0.69400,0.12500}%
\definecolor{mycolor4}{rgb}{0.49400,0.18400,0.55600}%

\begin{axis}[
/pgf/number format/1000 sep={},
width=4.521in,
height=3.566in,
at={(0.758in,0.481in)},
scale=0.4,
legend cell align={left},
legend style={
  fill opacity=0,
  draw opacity=1,
  text opacity=1,
  at={(1.1,0.97)},
  anchor=north west,
  draw=none
},
tick align=outside,
tick pos=left,
x grid style={black},
xlabel={Time (years)},
xmajorgrids,
xmin=2020, xmax=2052,
xtick style={color=black},
y grid style={black},
ylabel={$N_v$ (Mvehicles)},
ymajorgrids,
ymin=0, ymax=2,
ytick style={color=black}
]

\addplot [semithick, mycolor1, mark=*, mark size=.5, mark options={solid}]
table {%

2023 1.01517927904732
2024 1.06467183121248
2025 1.11519865799253
2026 1.1557135713313
2027 1.18471142318081
2028 1.20203520853632
2029 1.24291730044253
2030 1.2831198606974
2031 1.30174932808686
2032 1.32345233175703
2033 1.4154064917187
2034 1.43233137434477
2035 1.45625342319946
2036 1.48512134800008
2037 1.51724599382448
2038 1.55123387348811
2039 1.58592484972935
2040 1.62039571128848
2041 1.65384700243842
2042 1.68571720856938
2043 1.71560340930812
2044 1.74330315840776
2045 1.76868313381304
2046 1.79172890567911
2047 1.81252999716017
2048 1.83127810029943
2049 1.84824414494137
2050 1.81754743365601
};
\addlegendentry{EV (Opt.)}

\addplot [semithick, mycolor3, mark=*, mark size=.5, mark options={solid}]
table {%

2023 1.70315855846634
2024 1.57921548136421
2025 1.48032603311622
2026 1.39129158658586
2027 1.30967557788741
2028 1.23190130990965
2029 1.13478753837855
2030 1.04868262843489
2031 0.994418447912606
2032 0.947166516038613
2033 0.839136296413857
2034 0.814580815062154
2035 0.790137264871424
2036 0.766330126687121
2037 0.743304927651748
2038 0.720980456038655
2039 0.699177932177081
2040 0.677722776854636
2041 0.656449011236205
2042 0.635271536293914
2043 0.614156973973312
2044 0.593136030210035
2045 0.572251211560599
2046 0.551571730143368
2047 0.531183647157362
2048 0.511185095232278
2049 0.491676353888846
2050 0.518879603694141

};
\addlegendentry{ICEV (Opt.)}

\addplot [semithick, mycolor2, mark=*, mark size=.5, mark options={solid}]
table {%
2023 0.62580652378235
2024 0.733020415911391
2025 0.845749378103449
2026 0.946376442791526
2027 1.02582653254529
2028 1.08335491076081
2029 1.18873546521644
2030 1.29958796087198
2031 1.3905438807877
2032 1.48972475109184
2033 1.59262823621004
2034 1.6063756014762
2035 1.62690518534374
2036 1.65228975273203
2037 1.68088941469848
2038 1.71131195105727
2039 1.74237556984323
2040 1.77313400618992
2041 1.80276439570704
2042 1.83069824918656
2043 1.85653620962151
2044 1.88009327621429
2045 1.90125579787058
2046 1.92003479917163
2047 1.93654886404368
2048 1.9510212746768
2049 1.96375474999276
2050 1.9357186557451
};
\addlegendentry{EV (IC)}

\end{axis}

\end{tikzpicture}

%% file: figParc_aac/Stock5KOpt.tex
\definecolor{mycolor1}{rgb}{0.00000,0.44700,0.74100}%
\definecolor{mycolor2}{rgb}{0.85000,0.32500,0.09800}%
\definecolor{mycolor3}{rgb}{0.92900,0.69400,0.12500}%
\definecolor{mycolor4}{rgb}{0.49400,0.18400,0.55600}%

\begin{tikzpicture}

\definecolor{darkgray176}{RGB}{176,176,176}
\definecolor{lightgray204}{RGB}{204,204,204}
\definecolor{steelblue31119180}{RGB}{31,119,180}

\begin{axis}[
/pgf/number format/1000 sep={},
width=4.521in,
height=3.566in,
at={(0.758in,0.481in)},
scale=0.4,
legend cell align={left},
legend style={fill opacity=0.2, draw opacity=1, text opacity=1, draw=none,
at={(1.1,0.97)},
 anchor=north west},
tick align=outside,
tick pos=left,
x grid style={darkgray176},
xlabel={Time (years)},
xmajorgrids,
xmin=2020, xmax=2052,
xtick style={color=black},
y grid style={darkgray176},
ylabel={$S_v$ (Mvehicles)},
ymajorgrids,
ymin=0, ymax= 40,
ytick style={color=black}
]
\addplot [semithick, mycolor1, mark=*, mark size=.5, mark options={solid}]
table {%

2022 0.276883
2023 1.29135486794373
2024 2.35490365133818
2025 3.46841519676797
2026 4.62167958963476
2027 5.80261408076466
2028 6.99862304917957
2029 8.22592828929042
2030 9.47330873035929
2031 10.7086343782447
2032 11.9247898406515
2033 13.182615207678
2034 14.3988247635827
2035 15.5733888339143
2036 16.7055196332013
2037 17.7942620455312
2038 18.8388496202931
2039 19.8383784348856
2040 20.7922292117582
2041 21.7002550093439
2042 22.562786288021
2043 23.3805539520607
2044 24.154661719457
2045 24.8864443852979
2046 25.5774007359189
2047 26.2291441522406
2048 26.8433829691061
2049 27.4219125398571
2050 27.9204094245999
};
\addlegendentry{EV (Opt.)}

\addplot [semithick, mycolor3, mark=*, mark size=.5, mark options={solid}]
table {%

2022 35.519509
2023 34.7825329307659
2024 33.9617324389014
2025 33.083521296982
2026 32.158109419604
2027 31.1975795559451
2028 30.2145273269816
2029 29.1927309383006
2030 28.1434114606382
2031 27.0986988881433
2032 26.0657086131055
2033 24.9836005454302
2034 23.9356604008552
2035 22.9219178538336
2036 21.943160689837
2037 21.0003440247777
2038 20.0942343092667
2039 19.2257354659038
2040 18.3954667722428
2041 17.6035751698473
2042 16.8497301983353
2043 16.1332009534464
2044 15.4528837171812
2045 14.80744369445
2046 14.195382098919
2047 13.6150855496677
2048 13.0648457118528
2049 12.5428672321329
2050 12.0934735503995
};
\addlegendentry{ICEV (Opt.)}
\addplot [color=mycolor2, mark=*, mark options={solid, mycolor2}, mark size=.5]
  table[row sep=crcr]{%
2022	0.276883\\
2023	0.90198211267876\\
2024	1.63387948077211\\
2025	2.47794174631282\\
2026	3.42186901063984\\
2027	4.44391861113421\\
2028	5.52124728177362\\
2029	6.69710325793011\\
2030	7.97000666968697\\
2031	9.31230240582858\\
2032	10.7239404782539\\
2033	12.2002755811868\\
2034	13.6441437964197\\
2035	15.0546221112408\\
2036	16.4291370157325\\
2037	17.7640184007189\\
2038	19.0548684113387\\
2039	20.2967919404005\\
2040	21.4852832249128\\
2041	22.6165959604329\\
2042	23.6878800947605\\
2043	24.6972043125328\\
2044	25.6436017441817\\
2045	26.5269593070818\\
2046	27.3479553979408\\
2047	28.1079858433601\\
2048	28.8090953418799\\
2049	29.4539002713503\\
2050	30.0061194771407\\
};
\addlegendentry{EV (IC)}

\end{axis}

\end{tikzpicture}

%% file: figParc/Pareto.tex
\begin{tikzpicture}
\definecolor{mycolor1}{rgb}{0.00000,0.44700,0.74100}%
\begin{axis}[
/pgf/number format/1000 sep={},
width=4.521in,
height=3.566in,
at={(0.758in,0.481in)},
scale=0.4,
xlabel={$\mathcal{E}$(T) (Gt)},
ylabel={I(T) (G\euro)},
grid=both,
grid style={dashed, gray!30},
ymin=-50, 
ymax=700,
xmin=.7, 
xmax=1,
legend style={at={(0.95,0.95)}, anchor=north east},
]

\addplot[color=mycolor1,  mark=*, mark size=1, mark options={solid}] coordinates {
    (.98, 0)
    (.96, 26.5)
    (.91, 95.3)
    (.87, 185.2)
    (.82, 297)
    (.73, 622.3)
    };

\end{axis}
\end{tikzpicture}

%% file: main.bbl
\begin{thebibliography}{}

\bibitem [\protect \citeauthoryear {%
Ashina%
, Fujino%
, Masui%
, Ehara%
\BCBL {}\ \BBA {} Hibino%
}{%
Ashina%
\ \protect \BOthers {.}}{%
{\protect \APACyear {2012}}%
}]{%
Ashina}
\APACinsertmetastar {%
Ashina}%
\begin{APACrefauthors}%
Ashina, S.%
, Fujino, J.%
, Masui, T.%
, Ehara, T.%
\BCBL {}\ \BBA {} Hibino, G.%
\end{APACrefauthors}%
\unskip\
\newblock
\APACrefYearMonthDay{2012}{}{}.
\newblock
{\BBOQ}\APACrefatitle {A roadmap towards a low-carbon society in Japan using backcasting methodology: Feasible pathways for achieving an 80\% reduction in CO2 emissions by 2050} {A roadmap towards a low-carbon society in japan using backcasting methodology: Feasible pathways for achieving an 80\% reduction in co2 emissions by 2050}.{\BBCQ}
\newblock
\APACjournalVolNumPages{Energy Policy}{41}{}{584--598}.
\newblock
\begin{APACrefDOI} \doi{10.1016/j.enpol.2011.10.053} \end{APACrefDOI}
\PrintBackRefs{\CurrentBib}

\bibitem [\protect \citeauthoryear {%
Bass%
}{%
Bass%
}{%
{\protect \APACyear {1969}}%
}]{%
Bass}
\APACinsertmetastar {%
Bass}%
\begin{APACrefauthors}%
Bass, F\BPBI M.%
\end{APACrefauthors}%
\unskip\
\newblock
\APACrefYearMonthDay{1969}{}{}.
\newblock
{\BBOQ}\APACrefatitle {A New-Product Growth Model for Consumer Durables} {A new-product growth model for consumer durables}.{\BBCQ}
\newblock
\APACjournalVolNumPages{Management Science}{15}{1}{215--227}.
\PrintBackRefs{\CurrentBib}

\bibitem [\protect \citeauthoryear {%
Ben-Akiva%
\ \BBA {} Lerman%
}{%
Ben-Akiva%
\ \BBA {} Lerman%
}{%
{\protect \APACyear {2000}}%
}]{%
benakiva}
\APACinsertmetastar {%
benakiva}%
\begin{APACrefauthors}%
Ben-Akiva, L.%
\BCBT {}\ \BBA {} Lerman, S\BPBI R.%
\end{APACrefauthors}%
\unskip\
\newblock
\APACrefYear{2000}.
\newblock
\APACrefbtitle {Discrete Choice Analysis} {Discrete choice analysis}.
\newblock
\APACaddressPublisher{Cambridge, MA}{MIT Press}.
\PrintBackRefs{\CurrentBib}

\bibitem [\protect \citeauthoryear {%
Bibri%
\ \BBA {} Krogstie%
}{%
Bibri%
\ \BBA {} Krogstie%
}{%
{\protect \APACyear {2019}}%
}]{%
Bibri}
\APACinsertmetastar {%
Bibri}%
\begin{APACrefauthors}%
Bibri, S\BPBI E.%
\BCBT {}\ \BBA {} Krogstie, J.%
\end{APACrefauthors}%
\unskip\
\newblock
\APACrefYearMonthDay{2019}{}{}.
\newblock
{\BBOQ}\APACrefatitle {A scholarly backcasting approach to a novel model for smart sustainable cities of the future: strategic problem orientation} {A scholarly backcasting approach to a novel model for smart sustainable cities of the future: strategic problem orientation}.{\BBCQ}
\newblock
\APACjournalVolNumPages{City, Territory and Architecture}{6}{3}{1--27}.
\newblock
\begin{APACrefDOI} \doi{10.1186/s40410-019-0102-3} \end{APACrefDOI}
\PrintBackRefs{\CurrentBib}

\bibitem [\protect \citeauthoryear {%
Bouter%
\ \protect \BOthers {.}}{%
Bouter%
\ \protect \BOthers {.}}{%
{\protect \APACyear {2022}}%
}]{%
ADEME}
\APACinsertmetastar {%
ADEME}%
\begin{APACrefauthors}%
Bouter, A.%
\BCBT {}\ \BOthersPeriod {.}
\end{APACrefauthors}%
\unskip\
\newblock
\APACrefYearMonthDay{2022}{}{}.
\newblock
\APACrefbtitle {Etude énergétique, économique et environnementale du transport routier à horizon 2040 (E4T 2040)} {Etude énergétique, économique et environnementale du transport routier à horizon 2040 (e4t 2040)}\ \APACbVolEdTR{}{\BTR{}}.
\newblock
\APACaddressInstitution{}{IFPEN-ADEME}.
\PrintBackRefs{\CurrentBib}

\bibitem [\protect \citeauthoryear {%
De~Ceuster%
\ \protect \BOthers {.}}{%
De~Ceuster%
\ \protect \BOthers {.}}{%
{\protect \APACyear {2004}}%
}]{%
TREMOVE}
\APACinsertmetastar {%
TREMOVE}%
\begin{APACrefauthors}%
De~Ceuster, G.%
\BCBT {}\ \BOthersPeriod {.}
\end{APACrefauthors}%
\unskip\
\newblock
\APACrefYearMonthDay{2004}{}{}.
\newblock
\APACrefbtitle {TREMOVE 2.2 Model and Baseline Description} {Tremove 2.2 model and baseline description}\ \APACbVolEdTR{}{\BTR{}}.
\newblock
\APACaddressInstitution{}{European Commission}.
\PrintBackRefs{\CurrentBib}

\bibitem [\protect \citeauthoryear {%
{European Commission}%
}{%
{European Commission}%
}{%
{\protect \APACyear {2021}}%
}]{%
GreenDeal}
\APACinsertmetastar {%
GreenDeal}%
\begin{APACrefauthors}%
{European Commission}.%
\end{APACrefauthors}%
\unskip\
\newblock
\APACrefYearMonthDay{2021}{}{}.
\newblock
\APACrefbtitle {European Green Deal: Commission proposes transformation of EU economy and society to meet climate ambitions} {European green deal: Commission proposes transformation of eu economy and society to meet climate ambitions}\ \APACbVolEdTR{}{\BTR{}}.
\newblock
\APACrefnote{Press Release}
\PrintBackRefs{\CurrentBib}

\bibitem [\protect \citeauthoryear {%
Gomi%
, Ochi%
\BCBL {}\ \BBA {} Matsuoka%
}{%
Gomi%
\ \protect \BOthers {.}}{%
{\protect \APACyear {2011}}%
}]{%
Gomi}
\APACinsertmetastar {%
Gomi}%
\begin{APACrefauthors}%
Gomi, K.%
, Ochi, Y.%
\BCBL {}\ \BBA {} Matsuoka, Y.%
\end{APACrefauthors}%
\unskip\
\newblock
\APACrefYearMonthDay{2011}{}{}.
\newblock
{\BBOQ}\APACrefatitle {A systematic quantitative backcasting on low-carbon society policy in case of Kyoto city} {A systematic quantitative backcasting on low-carbon society policy in case of kyoto city}.{\BBCQ}
\newblock
\APACjournalVolNumPages{Technological Forecasting and Social Change}{78}{5}{852--871}.
\newblock
\begin{APACrefDOI} \doi{10.1016/j.techfore.2011.01.005} \end{APACrefDOI}
\PrintBackRefs{\CurrentBib}

\bibitem [\protect \citeauthoryear {%
{ITF}%
}{%
{ITF}%
}{%
{\protect \APACyear {2019}}%
}]{%
ITF}
\APACinsertmetastar {%
ITF}%
\begin{APACrefauthors}%
{ITF}.%
\end{APACrefauthors}%
\unskip\
\newblock
\APACrefYearMonthDay{2019}{}{}.
\newblock
\APACrefbtitle {Understanding Consumer Vehicle Choice: A New Car Fleet Model for France} {Understanding consumer vehicle choice: A new car fleet model for france}\ \APACbVolEdTR{}{\BTR{}}.
\newblock
\APACaddressInstitution{}{OECD Publishing, Paris}.
\PrintBackRefs{\CurrentBib}

\bibitem [\protect \citeauthoryear {%
Kenta%
\ \BBA {} Nakata%
}{%
Kenta%
\ \BBA {} Nakata%
}{%
{\protect \APACyear {2020}}%
}]{%
sato}
\APACinsertmetastar {%
sato}%
\begin{APACrefauthors}%
Kenta, S.%
\BCBT {}\ \BBA {} Nakata, T.%
\end{APACrefauthors}%
\unskip\
\newblock
\APACrefYearMonthDay{2020}{}{}.
\newblock
{\BBOQ}\APACrefatitle {Recoverability Analysis of Critical Materials from Electric Vehicle Lithium-Ion Batteries through a Dynamic Fleet-Based Approach for Japan} {Recoverability analysis of critical materials from electric vehicle lithium-ion batteries through a dynamic fleet-based approach for japan}.{\BBCQ}
\newblock
\APACjournalVolNumPages{Sustainability}{12}{1}{147}.
\PrintBackRefs{\CurrentBib}

\bibitem [\protect \citeauthoryear {%
McManus%
\ \BBA {} Senter%
}{%
McManus%
\ \BBA {} Senter%
}{%
{\protect \APACyear {2009}}%
}]{%
macmanus}
\APACinsertmetastar {%
macmanus}%
\begin{APACrefauthors}%
McManus, W.%
\BCBT {}\ \BBA {} Senter, R.%
\end{APACrefauthors}%
\unskip\
\newblock
\APACrefYearMonthDay{2009}{}{}.
\newblock
{\BBOQ}\APACrefatitle {Market Models for Predicting PHEV Adoption and Diffusion} {Market models for predicting phev adoption and diffusion}.{\BBCQ}
\newblock

\newblock
\APACrefnote{Final Report}
\PrintBackRefs{\CurrentBib}

\bibitem [\protect \citeauthoryear {%
Papazikou%
\ \protect \BOthers {.}}{%
Papazikou%
\ \protect \BOthers {.}}{%
{\protect \APACyear {2020}}%
}]{%
Levitate}
\APACinsertmetastar {%
Levitate}%
\begin{APACrefauthors}%
Papazikou, E.%
\BCBT {}\ \BOthersPeriod {.}
\end{APACrefauthors}%
\unskip\
\newblock
\APACrefYearMonthDay{2020}{}{}.
\newblock
\APACrefbtitle {Detailed list of sub-use cases, applicable forecasting methodologies and necessary output variables} {Detailed list of sub-use cases, applicable forecasting methodologies and necessary output variables}\ \APACbVolEdTR{}{\BTR{}}.
\newblock
\APACaddressInstitution{}{Deliverable D4.4 of the H2020 project LEVITATE}.
\PrintBackRefs{\CurrentBib}

\bibitem [\protect \citeauthoryear {%
Robinson%
}{%
Robinson%
}{%
{\protect \APACyear {1982}}%
}]{%
Robinson}
\APACinsertmetastar {%
Robinson}%
\begin{APACrefauthors}%
Robinson, J.%
\end{APACrefauthors}%
\unskip\
\newblock
\APACrefYearMonthDay{1982}{}{}.
\newblock
{\BBOQ}\APACrefatitle {Energy backcasting—A proposed method of policy analysis} {Energy backcasting—a proposed method of policy analysis}.{\BBCQ}
\newblock
\APACjournalVolNumPages{Energy Policy}{12}{4}{337--344}.
\newblock
\begin{APACrefDOI} \doi{10.1016/0301-4215(82)90048-9} \end{APACrefDOI}
\PrintBackRefs{\CurrentBib}

\bibitem [\protect \citeauthoryear {%
Struben%
\ \BBA {} Sterman%
}{%
Struben%
\ \BBA {} Sterman%
}{%
{\protect \APACyear {2008}}%
}]{%
sterman}
\APACinsertmetastar {%
sterman}%
\begin{APACrefauthors}%
Struben, J.%
\BCBT {}\ \BBA {} Sterman, J.%
\end{APACrefauthors}%
\unskip\
\newblock
\APACrefYearMonthDay{2008}{}{}.
\newblock
{\BBOQ}\APACrefatitle {Transition Challenges for Alternative Fuel Vehicle and Transportation Systems} {Transition challenges for alternative fuel vehicle and transportation systems}.{\BBCQ}
\newblock
\APACjournalVolNumPages{Environment and Planning B}{35}{6}{1070--1097}.
\PrintBackRefs{\CurrentBib}

\bibitem [\protect \citeauthoryear {%
Thorne%
, Aguilar~Lopez%
, Figenbaum%
, Fridstrøm%
\BCBL {}\ \BBA {} Müller%
}{%
Thorne%
\ \protect \BOthers {.}}{%
{\protect \APACyear {2021}}%
}]{%
thorne}
\APACinsertmetastar {%
thorne}%
\begin{APACrefauthors}%
Thorne, R.%
, Aguilar~Lopez, F.%
, Figenbaum, E.%
, Fridstrøm, L.%
\BCBL {}\ \BBA {} Müller, D\BPBI B.%
\end{APACrefauthors}%
\unskip\
\newblock
\APACrefYearMonthDay{2021}{}{}.
\newblock
{\BBOQ}\APACrefatitle {Estimating stocks and flows of electric passenger vehicle batteries in the Norwegian fleet from 2011 to 2030} {Estimating stocks and flows of electric passenger vehicle batteries in the norwegian fleet from 2011 to 2030}.{\BBCQ}
\newblock
\APACjournalVolNumPages{Journal of Industrial Ecology}{25}{6}{1377--1706}.
\PrintBackRefs{\CurrentBib}

\bibitem [\protect \citeauthoryear {%
Train%
}{%
Train%
}{%
{\protect \APACyear {2003}}%
}]{%
train}
\APACinsertmetastar {%
train}%
\begin{APACrefauthors}%
Train, K\BPBI E.%
\end{APACrefauthors}%
\unskip\
\newblock
\APACrefYear{2003}.
\newblock
\APACrefbtitle {Discrete Choice Methods with Simulation} {Discrete choice methods with simulation}\ (\PrintOrdinal{2nd}\ \BEd).
\newblock
\APACaddressPublisher{Cambridge}{Cambridge University Press}.
\PrintBackRefs{\CurrentBib}

\bibitem [\protect \citeauthoryear {%
Van~Grol%
\ \protect \BOthers {.}}{%
Van~Grol%
\ \protect \BOthers {.}}{%
{\protect \APACyear {2016}}%
}]{%
High-Tool}
\APACinsertmetastar {%
High-Tool}%
\begin{APACrefauthors}%
Van~Grol, R.%
\BCBT {}\ \BOthersPeriod {.}
\end{APACrefauthors}%
\unskip\
\newblock
\APACrefYearMonthDay{2016}{}{}.
\newblock
\APACrefbtitle {Elasticities and Equations of the HIGH-TOOL Model (Final Version)} {Elasticities and equations of the high-tool model (final version)}\ \APACbVolEdTR{}{\BTR{}}.
\newblock
\APACrefnote{HIGH-TOOL Deliverable D4.3, project cofunded by the European Commission under the 7th Framework Programme, Karlsruhe}
\PrintBackRefs{\CurrentBib}

\bibitem [\protect \citeauthoryear {%
Woody%
, Keoleian%
\BCBL {}\ \BBA {} Vaishnav%
}{%
Woody%
\ \protect \BOthers {.}}{%
{\protect \APACyear {2023}}%
}]{%
woody}
\APACinsertmetastar {%
woody}%
\begin{APACrefauthors}%
Woody, M.%
, Keoleian, G\BPBI A.%
\BCBL {}\ \BBA {} Vaishnav, P.%
\end{APACrefauthors}%
\unskip\
\newblock
\APACrefYearMonthDay{2023}{}{}.
\newblock
{\BBOQ}\APACrefatitle {Decarbonization potential of electrifying 50\% of U.S. light-duty vehicle sales by 2030} {Decarbonization potential of electrifying 50\% of u.s. light-duty vehicle sales by 2030}.{\BBCQ}
\newblock
\APACjournalVolNumPages{Nature Communications}{14}{}{7077}.
\PrintBackRefs{\CurrentBib}

\end{thebibliography}
